\documentclass[12pt]{article}
\usepackage{amsfonts,amssymb,amsmath}
\usepackage{theorem}
\newtheorem{definition}{Definition}[section]

\newtheorem{proposition}[definition]{Proposition}

\newtheorem{lemma}[definition]{Lemma}

\theorembodyfont{\rmfamily}

\newcommand{\CC}{{\mathbb C}}
\newcommand{\II}{{\mathbb I}}

\newcommand{\ZZ}{{\mathbb Z}}

\newcommand{\un}{{\mathbb I}}

\newcommand{\bi}{{\bar\imath}}
\newcommand{\bj}{{\bar\jmath}}

\newcommand{\sfrac}[2]{{\textstyle{\frac{#1}{#2}}}}
\newcommand{\half}{{\sfrac{1}{2}}}

\newcommand{\finproof}{{\hfill \rule{5pt}{5pt}}}

\newcommand{\qmbox}[1]{{\qquad\mbox{#1}\quad}}
\def\qmbox#1{\qquad\mbox{#1}\quad}

\def\tr{\mathop{\rm Tr}\nolimits}
\def\str{\mathop{\rm Str}\nolimits}
\def\diag{\mathop{\rm diag}\nolimits}
\newcommand{\non}{\nonumber}

\begin{document}
\pagestyle{empty}


\begin{center}

{\Large \textsf{Classification of reflection matrices related to
(super)~Yangians \\
and \\[3mm]
application to open spin chain models}}

\vspace{10mm}

{\large D. Arnaudon$^a$, J. Avan$^b$\footnote{On leave of absence from
  LPTHE, CNRS, UMR 7589,
  Universit{\'e}s Paris VI/VII.} , N.~Cramp\'e$^a$,
  A.~Doikou$^a$, L. Frappat$^{ac}$, {E}. Ragoucy$^a$}

\vspace{10mm}

\emph{$^a$ Laboratoire d'Annecy-le-Vieux de Physique Th{\'e}orique}

\emph{LAPTH, CNRS, UMR 5108, Universit{\'e} de Savoie}

\emph{B.P. 110, F-74941 Annecy-le-Vieux Cedex, France}

\vspace{7mm}

\emph{$^b$ Laboratoire de Physique Th{\'e}orique et Mod\'elisation}

\emph{Universit\'e de Cergy, 5 mail Gay-Lussac, Neuville-sur-Oise}

\emph{F-95031 Cergy-Pontoise Cedex}

\vspace{7mm}

\emph{$^c$ Member of Institut Universitaire de France}

\end{center}

\vfill
\vfill

\begin{abstract}
We present a classification of diagonal, antidiagonal and mixed reflection
matrices related to Yangian and super-Yangian $R$ matrices associated to
the infinite series $so(m)$, $sp(n)$ and $osp(m|n)$. We formulate the
analytical Bethe Ansatz resolution for the $so(m)$ and $sp(n)$ open spin
chains with boundary conditions described by the diagonal solutions.
\end{abstract}

\vfill
\vfill

\rightline{LAPTH-975/03}
\rightline{math.QA/0304150}
\rightline{March 2003}

\baselineskip=16pt


\newpage

\pagestyle{plain}
\setcounter{page}{1}

\section{Introduction}
\setcounter{equation}{0}

Quantum $R$ matrices, solutions of the Yang--Baxter equation, are
interpreted as diffusion amplitudes for two-body interactions of particle
type eigenstates in integrable two-dimensional field theories. The
Yang--Baxter equation then ensures consistent factorisability of the
three-body amplitudes in terms of two-body ones.

When considering integrable field theories with non trivial boundary
effects such as theories on a half-line, one needs to introduce a
new object describing reflection processes on the boundary. Integrability
is preserved provided that the two-body exchange $R$ matrix and the
one-body reflection $K$ matrix obey a quartic consistency condition
\cite{Che,Skl,GhZ,dVG}
\begin{equation}
  R_{12}(u-v)\, K_{1}(u)\, R_{12}(u+v)\, K_{2}(v) =
  K_{2}(v)\, R_{12}(u+v)\, K_{1}(u)\, R_{12}(u-v)\,.
  \label{eq:bybe}
\end{equation}
Let us remark that the reflection equation also appears when considering
generalisation of ZF algebras (i.e. Zamolodchikov--Faddeev algebras
\cite{ZF}) allowing the presence of a boundary. These generalised ZF
algebras ensure the total scattering matrix of the model to be unitary. In
this approach, the commutation relations of ZF generators implement an
operator $b(u)$ which obeys the reflection equation. The $K$ matrix would
then be a representation of this boundary operator $b$ \cite{Mintchev}.

Another point of view has recently been presented using a universal
construction of reflection algebras as twists of quantum algebras
\cite{DoMu,DoKuMu}. The exchange equation (\ref{eq:bybe}) would then
describe a trivial representation of the quantum generators, although it is
not clear to us that this interpretation is valid for spectral parameter
dependent solutions.
Furthermore, solutions with quantum degrees of freedom on the boundary
have also been obtained in \cite{BaKo}.

These reflection equations, or boundary Yang--Baxter equations, have
recently drawn attention and systematic ways of computing some solutions
for given $R$ matrices were derived e.g. in \cite{Ma,DeNe,Ne}, recovering and
extending previous results derived for instance in \cite{AR,MeNe,Li,Ga}. The
considered cases were related to $A_{1}^{(1)}$ trigonometric $R$ matrices
for all spins \cite{DeNe}, and to $A_{2}^{(2)}$ \cite{Kim} 
and $A_{n}^{(1)}$ vector
representations \cite{Ne}, yielding a wealth of new solutions. 

In this paper, we will consider the case of $R$ matrices
corresponding to vector representations of Yangians and super-Yangians,
i.e. rational $R$ matrices constructed in \cite{Dr,soya}. The complete
classification for $gl(n)$ in vector representation was given in
\cite{LAPTH844}: it appears that generically any solution is conjugated (by
a constant matrix) to a diagonal solution. \cite{McS} then derived a series
of solutions for $so(m)$ and $sp(n)$ based upon an ansatz proposed by
Cherednik \cite{Che}. We present here, for the Lie (super)algebra
series $so(m)$, $sp(n)$ and $osp(m|n)$, a classification of purely
diagonal, purely antidiagonal and mixed diagonal/antidiagonal solutions by
directly solving the boundary equation (\ref{eq:bybe})
for a one-dimensional boundary quantum space. Once the
diagonal case is exhausted (corresponding to ``flavour-preserving''
reflection matrices), the most natural extension to look for indeed 
consists of
reflection matrices which preserve \emph{pairs} of conjugate states
(according to (\ref{eq:defibar})). Particular solutions have already been
derived for the $so(m)$ algebra \cite{McS,Mori} and the $sp(n)$ algebra
\cite{McS}. They can all be identified with particular diagonal solutions
of our classification.

We must point out that there exists in fact another notion of reflection
equation, which is related to the definition of twisted Yangians
\cite{olsh,MNO}, and arises also in the theory of coideal algebras
described in \cite{Nou,MRS}. This equation reads
\begin{equation}
  R_{12}(u-v)\, K_{1}(u)\, R^{t_1}_{12}(-u-v)\, K_{2}(v) =
  K_{2}(v)\, R^{t_1}_{12}(-u-v)\, K_{1}(u)\, R_{12}(u-v)
  \label{eq:bybealt}
\end{equation}
where the transposition $t$ is defined below (see definition
\ref{def:trans}).
As pointed out in \cite{soya}, this equation is actually the same as
(\ref{eq:bybe}) when $R(u)$ is the $R$ matrix of Yangians of type
$so(n)$, $sp(n)$ or $osp(m|n)$.

One essential purpose in establishing such a classification of reflection
matrices is to use them in order to construct and eventually to solve spin
chain models with a variety of boundary conditions. We present here
explicit resolutions, using the analytical Bethe Ansatz method, of the
$so(n)$ and $sp(n)$ open spin chains with boundary conditions determined by
the diagonal solutions of the reflection equations. In particular, we
derive explicit formulae for the eigenvalues of the low-lying excitations
of hole type. We obtain in each case the full scattering matrix without
ambiguities (e.g. CDD factors) including bulk and boundary interactions.
This result has in fact a very general character, not limited to the
particular models considered here. This actually provides us with bulk
and boundary $S$ matrices relevant for integrable field theories with non
trivial boundary conditions, for which they represent ``universal'' $S$
matrices. 
We expect, in analogy with the bulk case, 
that such boundary $S$ matrices should correspond to $sp(n)$ or $so(m)$
Gross--Neveu model with certain boundaries. 
For this particular model, the study of boundary conditions
associated with the different $K$ matrices deserves further attention
but goes beyond the scope of the present paper. 

\section{Generalities}
\setcounter{equation}{0}

Let $gl(m|n)$ be the $\ZZ_{2}$-graded algebra of $(m+n) \!\times\!
(m+n)$ matrices $X_{ij}$ and $\theta_0=\pm 1$. The
$\ZZ_{2}$-gradation is defined by
$(-1)^{[i]} = \theta_0$ if $1 \le i \le m$ and $(-1)^{[i]} = -\theta_0$
if $m+1 \le i \le m+n$. In the following, we will always assume
that $n$ is even.
\begin{definition}
  For each index $i$, we introduce a sign $\theta_i$
  \begin{equation}
    \label{eq:deftheta}
    \theta_i = \begin{cases}
      +1 & \qmbox{for} 1\le i\le m+\frac{n}{2}  \cr
      -1 & \qmbox{for} m+\frac{n}{2}+1 \le i\le m+n
    \end{cases}
  \end{equation}
  and a conjugate index $\bar\imath$
  \begin{equation}
    \label{eq:defibar}
    \bar\imath =
    \begin{cases}
    m+1-i &\qmbox{for} 1\le i \le m \\
    2m+n+1-i &\qmbox{for} m+1\le i\le m+n
  \end{cases}
  \end{equation}
\end{definition}
In particular $\theta_{i} \theta_{\bar{\imath}} = \theta_{0} (-1)^{[i]} $.
\begin{definition}
  \label{def:trans}
  For $A = \sum_{ij} \; A^{ij} \;E_{ij}$, we
  define the transposition $t$ by
  \begin{equation}
    \label{eq:t}
    A^t = \sum_{ij} (-1)^{[i][j]+[j]} \theta_i \theta_j \; A^{ij} \,
    E_{\bar\jmath \bar\imath}
    = \sum_{ij}  \left(A^{ij}\right)^t \,
    E_{ij}
  \end{equation}
  It satisfies $(A^t)^t = A$ and, for $\CC$-valued bosonic matrices,
  $(AB)^t = B^t A^t$.
\end{definition}
As usual $E_{ij}$ denotes the elementary matrix with entry 1 in row $i$
and column $j$ and zero elsewhere.
\\
We shall use a graded tensor product, i.e. such that, for $a$, $b$,
$c$ and $d$ with definite gradings,
$(a\otimes b)(c\otimes d) = (-1)^{[b][c]} ac \otimes bd$.
\begin{definition}
  Let $P$ be the (super)permutation operator
  (i.e.  $X_{21}\equiv PX_{12}P$)
  \begin{equation}
    \label{eq:Pdef}
    P = \sum_{i,j=1}^{m+n} (-1)^{[j]} E_{ij} \otimes E_{ji}
  \end{equation}
  and let
  \begin{equation}
    \label{eq:Qdef}
        Q = \sum_{i,j=1}^{m+n} (-1)^{[i][j]} \theta_{i} \theta_{j}
        E_{\bar{\jmath}\bar{\imath}} \otimes E_{ji} = P^{t_{1}} \;.
  \end{equation}
  We define the $R$ matrix
  \begin{equation}
    \label{eq:RPQ}
    \displaystyle R(u) = \II + \frac{P}{u} - \frac{Q}{u+\kappa}
  \end{equation}
  with $2\kappa = (m-n-2)\theta_{0}$.
\end{definition}
The $R$ matrix (\ref{eq:RPQ}) satisfies the super Yang--Baxter equation
\begin{equation}
  \label{eq:YBE}
  R_{12}(u) \, R_{13}(u+v) \, R_{23}(v) = R_{23}(v) \, R_{13}(u+v)
  \, R_{12}(u)
\end{equation}
where the graded tensor product is understood.

The operators $P$ and $Q$ satisfy
\begin{equation}
  \label{eq:PQ}
  P^2 = \II, \qquad
  PQ = QP = \theta_{0} Q, \qquad \mbox{and} \qquad
  Q^2 = \theta_0 (m-n)Q \;.
\end{equation}
The $R$ matrix (\ref{eq:RPQ}) is known to yield the $osp(m|n)$ Yangian
\cite{soya}, and leads to the \emph{non-super} orthogonal (taking $n =
0$, $\theta_0=1$) and symplectic ($m=0$, $\theta_0=-1$) Yangians. For
obvious reasons, we will call labels as orthogonal (resp. symplectic),
indices $i$ which satisfy $1\leq i\leq m$ (resp. $m+1\leq i\leq m+n$).

\medskip

Although we will restrict ourselves to the diagonal, antidiagonal and mixed
cases for the reflection matrices $K$, the following lemma can be used to
get more general solutions.
\begin{lemma}
  \label{lemma:invar}
  Let $K(u)$ be a solution of the reflection equation (\ref{eq:bybe}) and
  $U$ such that $UU^t = 1$ be a (constant) matrix of the orthogonal,
  symplectic or orthosymplectic (super)group (depending upon the choice of
  $R$). Then $K^t(u)$, $U K(u) \, U^t$ and $U K^t(u) \, U^t$ are also
  solutions of (\ref{eq:bybe}).
\end{lemma}
The proof is straighforward, using the invariance of the $R$ matrix under
conjugation by $U$ and $R_{12}^{t_{1}t_{2}} = R_{12}$.

\section{Diagonal solutions of the reflection equation}
\setcounter{equation}{0}

In this section we consider invertible diagonal solutions for $K(u)$, i.e.
\begin{equation}
  \label{eq:diag}
  K(u) = \diag\Big(
  k_1(u),\cdots,k_m(u) \,;\,  k_{m+1}(u),\cdots,k_{m+n}(u) \Big)
\end{equation}
where the semicolon emphasises the splitting between orthogonal and
symplectic indices. Here the $k_i(u)$ are supposed to be analytic
$\CC$-functions of $u$, i.e. the boundary quantum space is
one-dimensional.  

\begin{proposition}
\label{prop:D}
  There are three families of generic diagonal solutions and two
  particular cases
  \begin{enumerate}
  \item[D1:]  Solutions of $sl(m+n)$ type, with one free parameter, for
  $m$ even
    \begin{equation}
      \label{eq:D1}
      k_i(u) = 1 \;,
      \qquad
      k_{\bar \imath}(u) = \frac{1+cu}{1-cu}\;,
      \qquad \forall\, i\in \{1,...,\frac{m}{2};m+1,...,m+\frac{n}{2}\}
    \end{equation}
    This solution has no extension to odd $m$. \\
        This solution is obviously invariant under the action of
        $SL(\frac{m}{2}|\frac{n}{2})$.
  \item[D2:] Solutions with three different values of $k_l(u)$,
  depending on one free parameter
    \begin{eqnarray}
      \label{eq:D2}
      &&
      k_{1}(u) = \frac{1+c_1 u}{1-c_1 u} \;,
      \qquad
      k_{m}(u) = \frac{1+c_{m}u}{1-c_{m}u}
      \nonumber\\[2mm]
      &&
      k_j(u) = 1    \qquad \forall\, j \neq 1,m
       \nonumber\\[2mm]
      &&
      \mbox{where} \qquad (\kappa-\theta_0) c_1 c_m + c_1 + c_m = 0
    \end{eqnarray}
    This solution does not hold for $m=0,1$. \\
    The moduli space of this solution is invariant under the action of
    $OSP(m-2|n) \times SO(2)$.
  \item[D3:] Solutions without any free continuous parameter (but with
  two integer parameters)
    \begin{eqnarray}
      \label{eq:D3}
      &&
      k_i(u) = k_{\bar\imath}(u) = 1 \qquad \forall\,
      i\in\{1,...,m_1;m+1,...,m+n_1\}
      \nonumber\\[2mm]
      &&
      k_{i}(u) = k_{\bar\imath}(u) = \frac{1+cu}{1-cu} \qquad \forall\,
      i\in\{m_1+1,...,m-m_1;m+n_1+1,...,m+n-n_1\}
      \nonumber\\[2mm]
      &&
      \mbox{where} \qquad c=\frac{2}{\kappa-\theta_0(2m_1-2n_1-1)}
    \end{eqnarray}
        This solution is invariant under the action of
        $OSP(2m_{1}|2n_{1}) \times OSP(m-2m_{1}|n-2n_{1})$.
  \item[D4:] In the particular case of $so(4)$, the solution takes the more
  general form:
  \begin{equation}
    \label{eq:D4}
    K(u) = \diag\bigg(1\;,\;\frac{1+c_2 u}{1-c_2 u}
    \;,\;\frac{1+c_3 u}{1-c_3 u}\;,\;
    \frac{1+c_2 u}{1-c_2 u}\;\frac{1+c_3 u}{1-c_3 u}
    \bigg)
  \end{equation}
  This solution contains the three generic solutions D1 ($c_2c_3=0$),
  D2 ($c_2+c_3=0$) and D3 ($c_2=c_3=\infty$).
  \item[D5:] In the particular case of $so(2)$, any function-valued
  diagonal matrix is solution.
  \end{enumerate}

  In each case D1--D4, infinite value of the parameter $c$ is allowed. It
  corresponds to the constant   value $\frac{1+cu}{1-cu}=-1$.

  These solutions are given up to a global normalisation and a relabelling
  of the indices, provided it preserves the
  orthogonal/symplectic splitting and all the sets of conjugate indices
  $\{i,\bar\imath\}$.

  Let us remind that the solutions for $so(m)$ and $sp(n)$ algebras are
  obtained by setting $n=0$ or $m=0$ respectively.
\end{proposition}

Solutions found in \cite{Mori} for $so(m)$ algebras ($m>2$) are
identified after a suitable change of basis with the sets D1 and D2.
Solutions found in \cite{McS} for $so(m)$ and $sp(n)$ algebras are
identified with the set D3 and limiting cases of D1.
The solution D3 with the special value $m_1=(m-1)/2$ (with $n=0$)
corresponds to 
the only non trivial rational limit of the set of trigonometric
solutions given in \cite{BFKZ}.

Remark that the generic solutions D1, D2, D3 are associated with symmetric
superspaces (as they were introduced in \cite{Zirn}, in a different context)
 based on $OSP(m|n)$, namely
$OSP(m|n)/SL(\frac{m}{2}|\frac{n}{2})$,
$OSP(m|n)/OSP({m}-{2}|{n})\times SO(2)$,
and $OSP(m|n)/OSP({m}-{2m_1}|{n}-2n_1)\times OSP({2m_1}|2n_1)$, respectively.
In the case of
Lie algebras, we recover the usual symmetric spaces (based on orthogonal 
and symplectic algebras).

\textbf{Proof:}
The projection of the reflection equation (\ref{eq:bybe})
on $E_{ij}\otimes E_{ji}$ for $i\neq j,\bj$ reads
\begin{equation}
  \label{eq:R2-1}
  \frac{1}{u+v}\; k_i(u) k_i(v) + \frac{1}{u-v} \; k_j(u)k_i(v)
  =   \frac{1}{u+v} \; k_j(u) k_j(v) + \frac{1}{u-v} \; k_j(v)k_i(u) \;,
\end{equation}
the solution of which is given by
\begin{equation}
  \label{eq:Fij}
  \frac{k_i(u)}{k_j(u)} = \frac{1+c_{ij}u}{1-c_{ij}u} \;, \qquad
  i\neq j,\bj \;.
\end{equation}
The cases $c_{ij}=0,\infty$ correspond to the constant ratios
$\frac{k_i(u)}{k_j(u)} = \pm 1$.
\\
For convenience, we will introduce $F_{ij}(u)\equiv \frac{k_i(u)}{k_j(u)}$
(well-defined since $K(u)$ is supposed to be invertible for generic
$u$).
When  $i\neq j,\bj$,  $F_{ij}$ is then given by (\ref{eq:Fij}) for
some $c_{ij}$ and obviously $c_{ji}=-c_{ij}$. When defined,
since $F_{ij}(u) F_{jk}(u) = F_{ik}(u)$,
the
parameters $c$ must also satisfy
\begin{equation}
  \label{eq:cocycle}
  \left\{
  \begin{array}{l}
    c_{ij} + c_{jk} + c_{ki} = 0  \\[2mm]
    c_{ij} c_{jk} c_{ki} = 0
  \end{array}
  \right.
  \qmbox{for}
  i\neq j,\bj \,, \quad j \neq k,\bar k \quad \mbox{and}
\quad k\neq i,\bi\,.
  \qquad  \qquad
\end{equation}
To any solution for $K(u)$, we associate a partition of
$\{1,\dots,m+n\}$ where the classes are defined by
$i\equiv j \ \Leftrightarrow \ k_i=k_j \ (\Leftrightarrow \ c_{ij}=0)$.
The constraints (\ref{eq:cocycle}) are
sufficient to conclude that the partition associated to any solution
has at most three different classes, except when $m+n=4$, where
it can in principle have four classes.
Note that in the case of $sl(m|n)$ where the constraints (\ref{eq:cocycle})
hold without any restriction on the
indices, $K(u)$ is always built with at most two different functions.
\\
More precisely, if $m+n\neq 4$, (\ref{eq:cocycle}) implies that the
partition of indices is constituted either of two subsets, or of three
subsets. In the last case, it can only be of the form
\begin{itemize}
\item[{(D2)}] \hspace{1cm}
  $\{i\}$, $\{\bi\}$ and $\{1,\dots,m+n\} \setminus\{i,\bi\}$
\end{itemize}
where $i$ is an orthogonal index, i.e. $i\in \{1,\dots,m\}$.
\\
Projecting the reflection equation (\ref{eq:bybe}) on
$E_{\bi \bj}\otimes E_{ij}$, one gets (for $i\neq j, \bj$) after
taking $u\rightarrow v$
\begin{eqnarray}
  \label{eq:R4etoile}
  &&
  \kappa \left( F_{\bi i}(u) - F_{\bj j}(u) \right)
  + \theta_0 \left( F_{ji}(u) - F_{ij}(u) \right)
  + (2u+\kappa) \left( F_{\bj i}(u) - F_{\bi j}(u) \right)
  \nonumber\\[3mm]
  &&\hspace{9cm}
  = \; \sum_{l} (-1)^{[l]} \left( F_{li}(u) - F_{lj}(u) \right)
\end{eqnarray}
This equation show in particular that
$F_{ij}=1 \Leftrightarrow F_{\bi \bj}=1$. A recursion allows us to write the
two-subset partitions as either
\begin{itemize}
\item[{(D1)}] \hspace{1cm}
  $I=\{i_1,i_2,\dots\}$ and $\overline I=\{\bi_1,\bi_2,\dots\}$,
  i.e. $i\in I \Leftrightarrow \bi\not\in I$
\item[or]
\item[{(D3)}] \hspace{1cm}
  $I=\{i_1,\bi_1,i_2,\bi_2\dots\}$ and $J=\{j_1,\bj_1,j_2,\bj_2,\dots\}$,
  i.e. $i\in I \Leftrightarrow \bi\in I$
\end{itemize}
Knowing the possible forms of the partition, one can evaluate the sum
in (\ref{eq:R4etoile}). This equation involves at most two different
functions $F_{mn}$. It provides the constraints on the parameters
$c_{mn}$.
A global check ensures that all the remaining projections do
not lead to new constraints.
\\
The cases $m+n=4$ are solved by direct computation. Only the
$so(4)$ case eventually exhibits a more general solution (D4).
\\
Finally, the case of $so(2)$ is special. In particular, $R(u)$ appears
to be diagonal. A direct computation shows that
all function-valued diagonal matrix is solution.
\finproof

Note that when $\kappa=0$ (as in the case of $so(2)$), the
spectral parameter can be rescaled before taking the limit
$\kappa\rightarrow 0$, and the corresponding $R$ matrix does not
involve the identity anymore \cite{KS}. This is the $R$ matrix used in
\cite{Mori} for $so(2)$.

\section{Antidiagonal solutions of the reflection equation}
\setcounter{equation}{0}
We now look for invertible solutions of the ``antidiagonal'' form
\begin{equation}
  \label{eq:antiK}
  K(u)=\sum\limits_{i=1}^{m+n} \ell_i(u) E_{i\bi} \;.
\end{equation}

\begin{proposition}
  Solutions of the reflection equation with antidiagonal terms exist only
  in the pure $so(2m)$ or $sp(2n)$ cases. They are constant and only
  restricted by the set of constraints
  \begin{equation}
    \label{eq:antidiag}
    \ell_i \ell_\bi = 1\;, \qquad \forall i.
  \end{equation}
  In the special case of $so(2)$, any function-valued antidiagonal
  matrix is solution.
\end{proposition}

The proof of this proposition will be given in the next section since
the antidiagonal solutions appear to be particular cases of the mixed
solutions discussed below.

\section{Mixed solutions of the reflection equation}
\setcounter{equation}{0}
We now look for invertible solutions with terms both in the diagonal and in
the antidiagonal part
\begin{equation}
  \label{eq:mixteK}
  K(u)=\sum\limits_{i=1}^{m+n} \left( k_i(u) E_{ii} \;+\; \ell_i(u)
  E_{i\bi} \right)
\end{equation}
with at least one non zero $\ell_i$.

\begin{lemma}
    The case of $so(2)$ being excluded, the solutions of the reflection
    equation of the form
    (\ref{eq:mixteK}) are all constant (up to a global normalisation
    function).

    In the special case of $so(2)$, the set of solutions of the
    reflection equation is the union of  function-valued diagonal
    matrices and  function-valued antidiagonal matrices.
\end{lemma}

\textbf{Proof:}
Projecting the reflection equations on the elementary matrices
$E_{pq}\otimes E_{rs}$, one gets
\begin{align}
  & \mbox {on $E_{ij} \otimes E_{j\bi}$ : } && \frac{(-1)^{[j]}}{u+v} \;
  \Big( k_{i}(u) \, \ell_{i}(v) + k_{\bi}(v) \, \ell_{i}(u)
  \Big) +
  \frac{(-1)^{[j]}}{u-v} \; \Big( k_{j}(u) \, \ell_{i}(v) - k_{j}(v) \,
  \ell_{i}(u) \Big) = 0
  \label{eq:ijjbi}
  \\
  & \mbox {on $E_{ij} \otimes E_{\bj i}$ : } && -\frac{(-1)^{[j]}}{u+v} \;
  \Big( k_{j}(u) \, \ell_{\bj}(v) + k_{\bj}(v) \, \ell_{\bj}(u)
  \Big) - \frac{(-1)^{[j]}}{u-v} \; \Big( k_{i}(u) \, \ell_{\bj}(v) -
  k_{i}(v) \, \ell_{\bj}(u) \Big) = 0
  \label{eq:ijbji}
  \\
  & \mbox {on $E_{ij} \otimes E_{ji}$ : } && \frac{(-1)^{[j]}}{u+v} \; \Big(
  k_{i}(u) \, k_{i}(v) + \ell_{i}(u) \, \ell_{\bi}(v) -
  k_{j}(u) \,
  k_{j}(v) - \ell_{\bj}(u) \, \ell_{j}(v) \Big) \nonumber \\
  &&& - \frac{(-1)^{[j]}}{u-v} \; \Big( k_{i}(u) \, k_{j}(v) - k_{j}(u) \,
  k_{i}(v) \Big) = 0
  \label{eq:ijji}
  \\
  & \mbox {on $E_{ij} \otimes E_{\bj\bi}$ : } && -\frac{(-1)^{[j]}}{u-v} \;
  \Big( \ell_{i}(u) \, \ell_{\bj}(v) - \ell_{i}(v) \,
  \ell_{\bj}(u) \Big) = 0
  \label{eq:ijbjbi}
  \\
  & \mbox {on $E_{ij} \otimes E_{ij}$ : } &&
  -\frac{(-1)^{[i][j]+[i]+[j]}}{u+v+\kappa} \; \theta_{i} \theta_{j} \theta_{0}
  \left\{ \Big( (-1)^{[i]} \ell_{i}(u) \, \ell_{\bj}(v) - (-1)^{[j]}
    \ell_{\bj}(u) \, \ell_{i}(v) \Big) \right. \nonumber \\
  &&& \hspace{110pt} \left. + \frac{1}{u-v} \; \Big( \ell_{i}(u) \,
    \ell_{\bj}(v) - \ell_{\bj}(u) \, \ell_{i}(v) \Big)
  \right\} = 0
  \label{eq:ijij}
  \\
  & \mbox {on $E_{ij} \otimes E_{\bi\bj}$ : } && (-1)^{[i][j]+[i]+[j]}
  \theta_{i}
  \theta_{j} \left\{ \frac{-1}{u+v+\kappa} \; \Big( k_{i}(u) \,
    k_{\bj}(v) - k_{j}(u) \, k_{\bi}(v) \Big) + \right. \nonumber \\
  &&& + \frac{1}{u-v+\kappa} \; \Big( k_{i}(u) \, k_{\bi}(v) - k_{j}(u)
  \, k_{\bj}(v) + (-1)^{[i]} \ell_{i}(u) \, \ell_{\bi}(v) - (-1)^{[j]}
  \ell_{\bj}(u) \, \ell_{j}(v) \Big) \nonumber \\
  &&& - \frac{1}{u+v+\kappa} \; \Big( \frac{\theta_{0}}{u-v} \; \big(
  k_{\bi}(u) \, k_{\bj}(v) - k_{\bj}(u) \, k_{\bi}(v) \big) + \frac{\str
    K(u)}{u-v+\kappa} \; \big( k_{\bi}(v) - k_{\bj}(v) \big) \Big)
  \nonumber \\
  &&& \left. + \frac{\theta_{0}}{(u+v)(u-v+\kappa)} \Big( k_{\bi}(u) \,
    k_{\bi}(v) - k_{\bj}(u) \, k_{\bj}(v) + \ell_{i}(u) \,
    \ell_{\bi}(v) - \ell_{\bj}(u) \, \ell_{j}(v) \Big) \right\}
  = 0\label{eq:ijbibj}
\end{align}
\begin{align}
  & \mbox {on $E_{ij} \otimes E_{\bi j}$ : } && (-1)^{[i][j]+[i]+[j]} \theta_{i}
  \theta_{j} \left\{ \frac{-1}{u+v+\kappa} \; \Big( k_{i}(u) \,
    \ell_{\bj}(v) - (-1)^{[j]} \theta_{0} \ell_{\bj}(u) \, k_{\bi}(v)
    \Big) + \right. \nonumber \\
  &&& - \frac{1}{u-v+\kappa} \; \Big( k_{j}(u) \, \ell_{\bj}(v) +
  (-1)^{[j]} \theta_{0} \ell_{\bj}(u) \, k_{j}(v) +
  \frac{\theta_{0}}{u+v}
  \; \big( k_{\bj}(u) \, \ell_{\bj}(v) + \ell_{\bj}(u) \, k_{j}(v)
  \big) \Big) \nonumber \\
  &&& \left. - \frac{1}{u+v+\kappa} \; \Big( \frac{\theta_{0}}{u-v} \;
    \big( k_{\bi}(u) \, \ell_{\bj}(v) - \ell_{\bj}(u) \, k_{\bi}(v)
    \big) - \frac{\str K(u)}{u-v+\kappa} \; \ell_{\bj}(v) \Big) \right\}
  = 0\label{eq:ijbij}
  \\
  & \mbox {on $E_{ij} \otimes E_{i\bj}$ : } && (-1)^{[i][j]+[i]+[j]} \theta_{i}
  \theta_{j} \left\{ \frac{1}{u+v+\kappa} \; \Big( k_{j}(u) \,
    \ell_{i}(v) - (-1)^{[i]} \theta_{0} \ell_{i}(u) \, k_{\bj}(v)
    \Big) + \right. \nonumber \\
  &&& + \frac{1}{u-v+\kappa} \; \Big( k_{i}(u) \, \ell_{i}(v) +
  (-1)^{[i]} \theta_{0} \ell_{i}(u) \, k_{i}(v) + \frac{\theta_{0}}{u+v}
  \; \big( k_{\bi}(u) \, \ell_{i}(v) + \ell_{i}(u) \, k_{i}(v)
  \big) \Big) \nonumber \\
  &&& \left. - \frac{1}{u+v+\kappa} \; \Big( \frac{\theta_{0}}{u-v} \;
    \big( \ell_{i}(u) \, k_{\bj}(v) - k_{\bj}(u) \, \ell_{i}(v)
    \big) + \frac{\str K(u)}{u-v+\kappa} \; \ell_{i}(v) \Big) \right\}
  = 0\label{eq:ijibj}
  \\
  & \mbox {on $1 \otimes E_{ii}$ : } && \frac{1}{u^2-v^2} \; \Big(
  \ell_{i}(u) \, \ell_{\bi}(v) - \ell_{\bi}(u) \,
  \ell_{i}(v)\Big) = 0
  \label{eq:1ii} \\
  &&& \nonumber \\
  & \mbox {on $1 \otimes E_{i\bi}$ : } && \frac{1}{u^2-v^2} \; \Big(
  k_{i}(u) \, \ell_{i}(v) + \ell_{i}(u) \, k_{\bi}(v) - k_{\bi}(u)
  \, \ell_{i}(v) - \ell_{i}(u) \, k_{i}(v) \Big)
  = 0\label{eq:1ibi}
\end{align}
Of course, when several indices $i,j,\bi,\bj$ coincide, the corresponding
equations merge into a single one. \\
Consider now a couple $(i,j)$ of indices such that $i,j,\bi,\bj$ are all
different. Then eq. (\ref{eq:ijbjbi}) implies
\begin{eqnarray}
  \ell_{i}(u) \, \ell_{\bj}(v) &=& \ell_{\bj}(u) \,
  \ell_{i}(v) \;,
  \label{eq:compa1}
\end{eqnarray}
the solution of which is given by $\ell_{i}(u) = \ell_{i}(0) \,
\ell(u)$ where $\ell(u)$ is an arbitrary function, which can be
factorised by a change of normalisation of the diagonal elements. Hence one
can restrict $K$ in (\ref{eq:mixteK}) to have only constant antidiagonal
elements without loss of generality.

\medskip

Suppose now that $\ell_{i} \ne 0$ for some $i$ (purely diagonal
solutions are not considered in this section). The general solution to the
system formed by eqs. (\ref{eq:ijjbi}) and (\ref{eq:1ibi}) is then given
by\footnote{Actually, we have restricted ourselves to meromorphic
functions on $\CC$.}
\begin{equation}
k_{i}(u) = \gamma u + k_i(0) \quad;\quad
k_{\bar{\imath}}(u) = \gamma u - k_i(0)
\quad\mbox{ and }\quad k_{j}(u) = - \gamma u + k_{j}(0)
\quad (j\neq i,\bi)
\end{equation}
where $\gamma$ is a constant. The inspection of the $E_{\bi\bi}
\otimes E_{i\bi}$ coefficient leads to $\gamma = 0$. \\
It follows that all diagonal terms are constant. Hence the mixed solutions
for the reflection matrix are necessarily constant. \\
The case of $so(2)$ is solved by direct calculation.
\finproof

We now specify the exact form of these solutions. We establish (the
case of $so(2)$ being excluded):

\medskip

\begin{proposition}
Mixed solutions of the reflection equation for $osp(m|n)$ exist only
when $m$ is even. They fall into two classes:
\begin{itemize}
\item[C1:] \underline{The $so(m)$ block is diagonal.}
\\[1.2ex]
The solutions are
parametrised by $n$ complex parameters. The matrix $K$ is given by:
\begin{equation}
\begin{array}{ll}
\mbox{for } i\in\{1,\ldots,\frac{m}{2}\} :
&\left\{\begin{array}{ll}
k_{i}=1 & k_{\bar\imath}=-1\\
\ell_{i}=0 & \ell_{\bar\imath}=0
\end{array} \right.\\[2.1ex]
\hspace{-2.1ex}\mbox{for } i\in\{m+1,\ldots,m+n\}:
& \left\{\begin{array}{l}
k_{i}=sin(\alpha_{i})\\
\ell_{i}=e^{\beta_{i}}\, cos(\alpha_{i})
\end{array}\right.\mbox{ with }
\begin{array}{l}
\alpha_{i}+\alpha_{\bar\imath}=0\\
\beta_{i}+\beta_{\bar\imath}=0
\end{array}
\end{array}
\end{equation}
\item[C2:] \underline{The $sp(n)$ block is diagonal.}
\\[1.2ex]
The solutions are parametrised by a couple of positive or null integers
$m_{1}\ge m_{2}$ such that
\begin{displaymath}
m_{1}+m_{2}\leq\frac{m}{2}-1\;;\quad\ m_{1}-m_{2} \leq \frac{n}{2}
\quad\mbox{and}\quad m_{1}-m_{2}\equiv \frac{n}{2}\ [mod\,2]
\end{displaymath}
and by $m-2(m_{1}+m_{2})$ complex parameters. Setting
$n_{1}=(n+2m_{1}-2m_{2})/4$, the matrix $K$ is given by:
\begin{equation}
\hspace{-2ex}so(m)\mbox{ part: }\left\{
\begin{array}{ll}
\mbox{for } i\in\{1,\ldots,m_{1}\}:
& \left[\begin{array}{l}
k_{i}=k_{\bar\imath}=1\\
\ell_{i}=\ell_{\bar\imath}=0
\end{array}\right.\\[2.1ex]
\mbox{for } i\in\{m_{1}+1,\ldots,m_{1}+m_{2}\}:
& \left[\begin{array}{l}
k_{i}=k_{\bar\imath}=-1\\
\ell_{i}=\ell_{\bar\imath}=0
\end{array}\right.\\[2.1ex]
\mbox{for } i\in\{m_{1}+m_{2}+1,\ldots,\frac{m}{2}\}:
& \left[\begin{array}{l}
k_{i}=k_{\bar\imath}=0\\
\ell_{\bar\imath}=\ell_{i}^{-1}\;,\ \ell_{i}\in\CC
\end{array}\right.
\end{array}\right.
\end{equation}
\begin{equation}
sp(n)\mbox{ part: }\left\{
\begin{array}{ll}
\mbox{for } i\in\{m+1,\ldots,m+n_{1}\}:
&\left[ \begin{array}{l}
k_{i}=k_{\bar\imath}=1\\
\ell_{i}=\ell_{\bar\imath}=0
\end{array} \right.\\[2.1ex]
\mbox{for } i\in\{m+n_{1}+1,\ldots,m+\frac{n}{2}\}:
&\left[\begin{array}{l}
k_{i}=k_{\bar\imath}=-1\\
\ell_{i}=\ell_{\bar\imath}=0
\end{array} \right.
\end{array}\right.
\end{equation}
\end{itemize}
For $so(m)$ (resp. $sp(n)$), the mixed solutions are of the form C2 with
$n=0$ and $m_{2}=m_{1}$ (resp. C1 with $m=0$).
\\
Note that this classification is given up to a global normalisation
function and up to a relabelling of the indices, as noticed in
proposition \ref{prop:D}.
\end{proposition}
\textbf{Proof:}
We have to consider only constant reflection matrices $K$. Then, extracting
the residues at the poles in $u$, eqs. (\ref{eq:ijjbi}) to (\ref{eq:1ibi})
reduce to (for $i \ne \bi,j,\bj$ and $j \ne \bj$)
\begin{eqnarray}
  && \ell_{i} \, (k_{i} + k_{\bi}) = 0
  \label{eq:c1} \\
  && k_{i}^2 - k_{j}^2 + \ell_{i} \ell_{\bi} - \ell_{j}
  \ell_{\bj} = 0
  \label{eq:c2} \\
  && k_{i} k_{\bi} - k_{j} k_{\bj} + \theta_{0} ((-1)^{[i]} \ell_{i}
  \ell_{\bi} - (-1)^{[j]} \ell_{j} \ell_{\bj}) = 0
  \label{eq:c3} \\
  && ((-1)^{[j]} - (-1)^{[i]}) \ell_{i} \ell_{j} = 0
  \label{eq:c4} \\
  && k_{i} k_{\bj} - k_{j} k_{\bi} = 0 \label{eq:c5} \\
  && (k_{\bi} - (-1)^{[j]} \theta_{0} k_{i}) \; \ell_{j} = 0
  \label{eq:c6} \\
  && (1 + (-1)^{[i]} \theta_{0}) \; k_{i} \ell_{i} = 0
  \label{eq:c7} \\
  && (1 + (-1)^{[i]} \theta_{0}) \; k_{i} \ell_{\bi} = 0
  \label{eq:c8} \\
  && (k_{i} - k_{j})\, \str K  = 0 \label{eq:c9} \\
  && \ell_{j} \, \str K = 0
  \label{eq:c10}
\end{eqnarray}

In the case of $so(m)$ algebras, eq. (\ref{eq:c6}) implies $k_{i} =
k_{\bi}$ for all $1 \le i \le m$, from which (\ref{eq:c5}) follows. Eq.
(\ref{eq:c1}) is a consequence of (\ref{eq:c7}) and (\ref{eq:c8}), eq.
(\ref{eq:c9}) follows from (\ref{eq:c10}) and (\ref{eq:c3}) from
(\ref{eq:c2}). Finally the mixed solutions in the $so(m)$ case are
characterised by the following equations:
\begin{eqnarray}
  && k_{i} = k_{\bi} \;; \qquad k_{i} \ell_{i} = 0 \;; \qquad \tr K
  = 0 \label{eq:solso1} \\
  && k_{i}^2 + \ell_{i} \ell_{\bi} = k_{j}^2 + \ell_{j}
  \ell_{\bj} \qquad \mbox{for all } i,j
  \label{eq:solso2}
\end{eqnarray}
These constraints exclude the existence of mixed or antidiagonal solutions
with odd $m$, since a matrix cannot be traceless if its diagonal has an odd
number of non vanishing elements with equal squares.

In the case of $sp(n)$ algebras, eqs. (\ref{eq:c7}) and (\ref{eq:c8})
vanish, while eq. (\ref{eq:c6}) implies $k_{i} = -k_{\bi}$ for all $1 \le i
\le n$, from which (\ref{eq:c1}), (\ref{eq:c5}), (\ref{eq:c9}) and
(\ref{eq:c10}) follow. Moreover, eqs. (\ref{eq:c2}) and (\ref{eq:c3}) are
equal. Finally the mixed solutions in the $sp(n)$ case are
characterised by the following equations:
\begin{eqnarray}
  && k_{i} = -k_{\bi} \;; \qquad \tr K = 0 \label{eq:solsp1} \\
  && k_{i}^2 + \ell_{i} \ell_{\bi} = k_{j}^2 + \ell_{j}
  \ell_{\bj} \qquad \mbox{for all } i,j
  \label{eq:solsp2}
\end{eqnarray}

In the case of $osp(1|n)$, one has to add the following equations (with $2
\le i,\bi \le 1+n$)
\begin{eqnarray}
  && k_{1}^2 - k_{i}^2 - \ell_{i} \ell_{\bi} = 0
  \label{eq:c11} \\
  && k_{1}^2 - k_{i} k_{\bi} - \ell_{i} \ell_{\bi} = 0
  \label{eq:c12}
\end{eqnarray}
It follows from (\ref{eq:c11})--(\ref{eq:c12}) and (\ref{eq:c8}) that
$k_{i}^2 = 0$ for all $i=2,\ldots,1+n$. Since $\ell_{i} \, \str K=
0$, one has also $k_{1} = 0$. Therefore, there does not exist any
mixed or antidiagonal (invertible) solution for $osp(1|n)$.

\medskip

In the case of $osp(m|n)$, ($m>1$), eq. (\ref{eq:c4}) shows that at least
one antidiagonal part (orthogonal or symplectic) is zero: $\ell_{i} = 0$
for $i \in [1,m]$ or $\ell_{i} = 0$ for $i \in [m+1,m+n]$. This means that
no pure antidiagonal solution exists for $osp(m|n)$ superalgebras. We
thus have to consider two cases. \\
$\bullet$ Case 1 : the $so(m)$ block is diagonal \\
Eq. (\ref{eq:c6}) with index $j \in [m+1,m+n]$ implies $k_{i} = -k_{\bi}$
for all $i$, which excludes the case of $osp(m|n)$ with odd $m$. The
remaining equations then lead to
\begin{eqnarray}
  && k_{i}^2 + \ell_{i} \ell_{\bi} = k_{j}^2 + \ell_{j}
     \ell_{\bj} \qquad \mbox{for } i,j \in [1,m+n]
     \label{eq:solosp2}
     \\
  && k_j = - k_\bj \qquad \mbox{for } j \in [1,m+n]
     \label{eq:solosp1}
\end{eqnarray}
\\
$\bullet$ Case 2 : the $sp(n)$ block is diagonal \\
Eq. (\ref{eq:c6}) with index $j \in [1,m]$ implies $k_{i} = k_{\bi}$ for
all $i$. The remaining equations lead then to
\begin{eqnarray}
  && k_{i}^2 + \ell_{i} \ell_{\bi} = k_{j}^2 + \ell_{j}
     \ell_{\bj} \qquad \mbox{for } i,j \in [1,m+n]
     \label{eq:solosp3}
     \\
  && k_j = k_\bj \qquad \mbox{for } j \in [1,m+n]
     \\
  && k_i \ell_i = k_i \ell_\bi = 0 \qquad \mbox{for } i \in [1,m]
     \\
  && \str K = \Big( \sum_{i=1}^{m} - \sum_{i=m+1}^{m+n} \Big) k_{i} = 0
     \label{eq:solosp4}
\end{eqnarray}
Again, these constraints exclude the existence of mixed or antidiagonal
solutions with odd $m$, since one cannot have a traceless matrix if its
diagonal has an odd number of non vanishing elements with equal squares.

\finproof
\medskip

For example, one finds the two following solutions for $osp(4|2)$:
\begin{displaymath}
  \left(
     \begin{array}{cccc|cc}
      {1} & & & & & \\
      & {1} & & & & \\
      & &  {-1} & & & \\
      & & & -{1} & & \\
      \hline
      \Big.
      & & & & k_{5} & \ell_5 \\[1mm]
      & & & & \ell_6 & - k_{5}
     \end{array}
   \right)
   \qquad \mbox{and} \qquad
   \left(
     \begin{array}{cccc|cc}
       1 & & & 0 & & \\
       & 0 & \ell_2 & & & \\
       & \ell_2^{-1} & 0 & & & \\
       0 & & & 1 & & \\
       \hline
       & & & & 1 & \\
       & & & & & 1
     \end{array}
   \right)
\end{displaymath}
where $k_{5}^2 + \ell_{5} \ell_{6} = 1$.
For $osp(2|4)$ the two solutions take the form
\begin{displaymath}
  \left(
    \begin{array}{cc|cccc}
      1 & & & & & \\
      & -1 & & & & \\
      \hline
      & & k_{3} & & & \ell_{3} \\
      & & & k_{4} & \ell_{4} & \\
      & & & \ell_{5} & -k_{4} & \\
      & & \ell_{6} & & & -k_{3}
    \end{array}
  \right)
  \qquad \mbox{and} \qquad
  \left(
    \begin{array}{cc|cccc}
      0 & \ell_1 & & & & \\[1mm]
      \ell_1^{-1} & 0 & & & & \\
      \hline
      & & 1 & & & \\
      & & & -1 & & \\
      & & & & -1 & \\
      & & & & & 1
    \end{array}
  \right)
\end{displaymath}
where $k_{3}^2 + \ell_{3} \ell_{6} = 1$ and
$k_{4}^2 + \ell_{4} \ell_{5} = 1$.

\section{Analytical Bethe Ansatz for the $\mathbf{so(n)}$ and
$\mathbf{sp(n)}$ open spin chains}
\setcounter{equation}{0}

The main aim of this section is the derivation of the Bethe Ansatz
equations for the $so(n)$ and $sp(n)$ $N$-site open spin chains with non
trivial reflection conditions by means of the analytical Bethe Ansatz
method (see e.g. \cite{reshe,mnanal,KunSuz,AMN1,AMN2,doikou}). To
construct the open chain 
transfer matrix we need to introduce the $R$ matrix \cite{baxter,korepin},
which is a solution of the Yang--Baxter equation (\ref{eq:YBE}).
Normalisation of $R$ matrix will be modified in order to connect it with
the physical $XXX$ type Hamiltonians considered hereafter; in addition, we
set $\lambda = iu$. We focus on the $so(n)$ or $sp(n)$ invariant $R$ matrix
given by:
\begin{eqnarray}
R(\lambda) =
\lambda(\lambda+i\kappa)\un +i(\lambda +i\kappa) P -i\lambda Q.
\label{r}
\end{eqnarray}
The $R$ matrix (\ref{r}) satisfies also crossing and unitarity, namely
\begin{eqnarray}
R_{12}(\lambda) R_{12}(-\lambda) = (\lambda^2+\kappa^2)(\lambda^2+1)\, \un,
~~R_{12}(\lambda) = R_{12}^{t_{1}}(-\lambda -i\kappa).
\label{cu}
\end{eqnarray}
where $^t$ is the transposition defined as in (\ref{eq:t}). \\
The open chain transfer matrix is defined by \cite{Skl}
\begin{eqnarray}
t(\lambda) = \tr_{0}  K_{0}^{+}(\lambda)\
T_{0}(\lambda)\ K^{-}_{0}(\lambda)\ \hat T_{ 0}(\lambda)\,,
\label{transfer1}
\end{eqnarray}
where $K^{-}_{0}(\lambda)$ is any solution of eq. (\ref{eq:bybe}) and
$K^{+}_{0}(\lambda)$ is a solution of a closely related reflection equation
defined to be:
\begin{equation}
  R_{12}(v-u)\, K_{1}^{t_{1}}(u)\, R_{12}(-u-v-2i\kappa)\, K_{2}^{t_{2}}(v) =
  K_{2}^{t_{2}}(v)\, R_{12}(-u-v-2i\kappa)\, K_{1}^{t_{1}}(u)\, R_{12}(v-u)
  \label{eq:bybeplus}
\end{equation}
$\tr_{0}$ denotes trace over the auxiliary space 0, and
\begin{eqnarray}
  T_{0}(\lambda) = R_{0N}(\lambda)  R_{0\,N-1}(\lambda) \cdots
  R_{0 2}(\lambda) R_{01}(\lambda)\,,~~
  \hat T_{0}(\lambda) =
  R_{10}(\lambda)  R_{20}(\lambda) \cdots
  R_{N-1\;0}(\lambda) R_{N0}(\lambda)\,
  \label{hatmonodromy}
\end{eqnarray}
It is clear that any solution $K^{-}(\lambda,\xi^{-})$ of (\ref{eq:bybe})
where $\xi^{-}$ are arbitrary boundary parameters, give rise to a solution
$K^{+}(\lambda)$ of (\ref{eq:bybeplus}) defined by $K^{+}(\lambda) =
K^{-}(-\lambda -i\kappa)^{t}$, also depending on
arbitrary boundary parameters $\xi^{+}=\xi^{-}$.

To determine the eigenvalues of the transfer matrix and the corresponding
Bethe Ansatz equations, we employ the analytical Bethe Ansatz method
\cite{reshe,mnanal}. Namely, we impose certain constraints on the
eigenvalues by exploiting the crossing symmetry of the model, the symmetry
of the transfer matrix, the analyticity of the eigenvalues, and the fusion
procedure for open spin chains. This then allows us to determine the
eigenvalues by solving a set of coupled non-linear consistency equations or
Bethe Ansatz equations. \\
We first focus on the case with trivial boundaries, namely $K^{-}(\lambda)
=K^{+}(\lambda)=1$. We will in a second step derive the Bethe Ansatz
equations for all diagonal solutions found in the previous sections, namely
D1, D2, D3, D4 and solve them in the thermodynamical limit.

To obtain the necessary constraints, we recall that the fusion procedure
for the open spin chain \cite{mnfusion,Zhou,doikou} yields the fused transfer
matrix
\begin{eqnarray}
\tilde t(\lambda) = \zeta(2\lambda+2i\kappa)\ t(\lambda)\ t(\lambda +
i\kappa) - (\zeta(\lambda+i\kappa))^{2N} q(2\lambda +i\kappa)q(-2 \lambda -
3i\kappa)\,,
\label{fusion}
\end{eqnarray}
where we define
\begin{eqnarray}
\zeta(\lambda) &\!\!=\!\!&
(\lambda+i\kappa)(\lambda+i)(\lambda-i\kappa)(\lambda-i) \,, \\
q(\lambda) &\!\!=\!\!&
\begin{cases}
(\lambda - i)(\lambda - i\kappa) & ~~ \mbox{for} ~~so(n) \\
(\lambda + i)(-\lambda + i\kappa) & ~~\mbox{for} ~~sp(n)
\end{cases}
\end{eqnarray}
In addition, from the crossing symmetry of the $R$ matrix (\ref{cu}) it
follows that (when $K^{-}=K^{+} =1$):
\begin{eqnarray}
t(\lambda) =t(-\lambda -i\kappa).
\label{cross}
\end{eqnarray}
The transfer matrix with $K^{-}=K^{+} =1$ is obviously $so(n)$ (resp.
$sp(n)$) invariant since the corresponding $R$ matrix (\ref{r}) is $so(n)$
(resp. $sp(n)$) invariant. The symmetry of the transfer matrix makes the
computation of its asymptotic behaviour ($\lambda \to \infty$) a relatively
easy task. Finally, to implement the analyticity of the eigenvalues, we
require that all poles must vanish. These constraints shall uniquely fix
the eigenvalues. This is the basic outline of the analytical Bethe Ansatz
method.

\subsection{Eigenvalue of the pseudo-vacuum}

In order to compute the general eigenvalues we need to first define a
reference state or ``pseudo-vacuum''. After finding its pseudo-energy
eigenvalue, we will be able, with the help of the above discussed
constraints, to derive the general eigenvalues of the transfer matrix. The
pseudo-vacuum, which is an exact eigenstate of the transfer matrix, is the
state with all ``spins'' up, i.e.
\begin{eqnarray}
\vert \omega_{+} \rangle = \bigotimes_{i=1}^{N} \vert + \rangle _{i}
\label{pseudo}
\end{eqnarray}
where $\vert + \rangle$ is the $n$-dimensional column vector
\begin{eqnarray}
\vert + \rangle = \left (
\begin{array}{c}
  1 \\
  0 \\
  \vdots  \\
  0 \\
\end{array}
\right)\,.
\label{col}
\end{eqnarray}
The action of the $R$ matrix on the state $\vert + \rangle$ from the left
gives rise to upper triangular matrices, whereas the action from the right
on $\langle + \vert = \vert + \rangle^\dagger$ gives lower triangular
matrices. It is obvious then that the corresponding action of the $\hat T$
(resp. $T$) on $\vert \omega_{+} \rangle$ (resp. $\langle \omega_{+}\vert$)
will also give rise to upper (resp. lower) triangular matrices. It is
relatively easy, after some tedious algebra (see \cite{doikou} for a
detailed example), to determine the action of the transfer matrix,
$t(\lambda) \vert \omega_{+} \rangle = \Lambda^{0}(\lambda) \vert
\omega_{+} \rangle$, where $\Lambda^{0}(\lambda)$ is given by the following
expression
\begin{eqnarray}
  \Lambda^{0}(\lambda) = a(\lambda)^{2N} g_{0}(\lambda)
  +b(\lambda)^{2N}\sum_{l=1}^{n-2}g_{l}(\lambda)+ c(\lambda)^{2N}
  g_{n-1}(\lambda)
  \label{eigen0}
\end{eqnarray}
with
\begin{eqnarray}
  a(\lambda) =(\lambda+i)(\lambda+i\kappa), ~~b(\lambda) =\lambda
  (\lambda+i\kappa), ~~c(\lambda) =\lambda (\lambda+i\kappa-i)
  =a(-\lambda-i\kappa)\quad
  \label{numbers}
\end{eqnarray}
\begin{eqnarray}
  g_{0}(\lambda) = {(\lambda+ {i\kappa \over 2} \pm {i \over
      2})(\lambda+i\kappa) \over(\lambda+ {i \kappa \over 2})(\lambda+{i \over
      2})}, ~~g_{n-1}(\lambda) ={(\lambda+ {i\kappa \over 2} \mp {i \over 2})
    \lambda \over(\lambda+ {i \kappa \over 2})(\lambda+ i\kappa -{i \over
      2})}=g_0(-\lambda-i\kappa)\quad \,.
  \label{g1}
\end{eqnarray}
In $g_{0}$ and $g_{n-1}$ the upper sign corresponds to $so(n)$ and the
lower sign to $sp(n)$. Moreover (when $n=2k$),
\begin{eqnarray}
g_{l}(\lambda) &=& {\lambda (\lambda+ {i\kappa \over 2} \pm {i \over
2})(\lambda+i\kappa) \over (\lambda+ {i \kappa \over 2})(\lambda+{i l \over
2})(\lambda+{i (l +1)\over 2})}, ~~ 0<l<k, \non\\
g_{l}(\lambda) &=& g_{n-l-1}(-\lambda -i\kappa), ~~ k \le l < n-1
\label{g2}
\end{eqnarray}
This expression is valid for $so(2k)$ with the upper sign 
and for $sp(2k)$ with the lower sign. For the $so(2k+1)$ chain the
expressions for $g_{l}$, $l\neq k$, are the same as in the $so(2k)$
chain, except for 
\begin{eqnarray}
g_{k}(\lambda) &=& {\lambda (\lambda+i\kappa) \over (\lambda+{i k \over
2})(\lambda+{i (k -1)\over 2})}= g_{k}(-\lambda -i\kappa).
\label{g3}
\end{eqnarray}

\subsection{Dressing functions}

Now that we have the expression for the pseudo-vacuum eigenvalue, in
accordance with the general analytical Bethe Ansatz procedure, we make the
following assumption for the structure of the general eigenvalues:
\begin{eqnarray}
\Lambda(\lambda) =
a(\lambda)^{2N} g_{0}(\lambda)A_{0}(\lambda)
+b(\lambda)^{2N}\sum_{i=1}^{n-2}g_{i}(\lambda)A_{i}(\lambda)+
c(\lambda)^{2N} g_{n-1}(\lambda)A_{n-1}(\lambda)
\label{eigen}
\end{eqnarray}
where $A_{i}(\lambda)$, the so-called ``dressing functions'', will now be
determined. \\
We immediately get from the crossing symmetry of the transfer matrix
(\ref{cross}):
\begin{eqnarray}
A_{l}(\lambda) = A_{n-l-1}(-\lambda
-i\kappa)\qquad 0 \le l \le n-1 \,.
\label{11}
\end{eqnarray}
Moreover, from the fusion relation (\ref{fusion}), we obtain the following
identity by a comparison of the forms (\ref{eigen}) for the initial and
fused auxiliary spaces:
\begin{eqnarray}
A_{0}(\lambda+i\kappa)A_{n-1}(\lambda) =1 \,.
\label{22}
\end{eqnarray}
Gathering the above two equations (\ref{11}), (\ref{22}) we conclude
\begin{eqnarray}
A_{0}(\lambda)A_{0}(-\lambda)=1 \,.
\label{33}
\end{eqnarray}
Additional constraints are obtained on the dressing functions from
analyticity properties. Studying carefully the common poles of successive
$g_{l}$'s, we deduce from the form of the $g_{i}$ functions
(\ref{g1})--(\ref{g3}) that $g_{l}$ and $g_{l-1}$ have common poles at
$\lambda =-{il\over 2}$, therefore from analyticity requirements
\begin{eqnarray}
A_{l}(-{il \over 2})=A_{l-1}(-{il \over 2}),
~~~l =1,\ldots,k-1 \,. \label{anal1}
\end{eqnarray}
The last relation is valid for $so(n)$ and $sp(n)$ as well. However, there
is one extra constraint specific to $so(2k+1)$, namely
\begin{eqnarray}
  A_{k}(-{ik \over 2})=A_{k-1}(-{ik \over 2}). 
  \label{anal2}
\end{eqnarray}
Having deduced the necessary constraints for the dressing functions, we
now determine them explicitly. The dressing functions $A_{l}$ are
essentially characterised by a set of parameters $\{\lambda_{j}^{(l)} \;|\;
j=1, \ldots, M^{(l)}\}$, where the integer numbers $M^{(l)}$ are related to
the eigenvalues of diagonal generators of $so(n)$, $sp(n)$. Defining these
generators as:
\begin{eqnarray}
S^{(l)}
= \sum_{i=1}^{N} s_{i}^{(l)}, ~~s^{(l)} =E_{ll} - E_{\bar l \bar l}. 
\label{qn}
\end{eqnarray}
The precise identification of $M$ follows from the symmetry of the
transfer matrix (see also \cite{reshe}), in particular
\begin{eqnarray}
S^{(l)} = M^{(l-1)} - M^{(l)}, ~~l=1,\ldots,k-2
\label{quant}
\end{eqnarray}
which is valid for $so(2k+1)$, $so(2k)$ and $sp(2k)$.
The two remaining quantum numbers are given by
\begin{eqnarray}
  \label{eq:quant2}
  &&\mbox{For } so(2k+1): \qquad S^{(k-1)} = M^{(k-2)} - M^{(k-1)} \;; \quad
  S^{(k)} = M^{(k-1)} - M^{(k)} \;;  \non\\
  &&\mbox{For } so(2k): \qquad S^{(k-1)} = M^{(k-2)} - M^{(+)} -
  M^{(-)}\;;\quad  S^{(k)} = M^{(+)}-M^{(-)} \;; \non\\
  &&\mbox{For } sp(2k): \qquad S^{(k-1)} = M^{(k-2)} - M^{(k-1)} \;; \quad
  S^{(k)} = M^{(k-1)} - 2M^{(k)} \;.  
\end{eqnarray}
We have also $M^{(0)} =N$ (valid for both $sp(n)$ and $so(n)$
algebras). \\
Let us point out that, away from the boundaries, the behavior of the chain
is considered to be as in the bulk. Therefore, the above quantum numbers
describe accurately enough the states of the system. The dressing functions
are given by the following expressions:

\subsection*{A. $\bf so(2k+1)$}
\begin{eqnarray}
A_{0}(\lambda) &=&  \prod_{j=1}^{M^{(1)}}
{\lambda+\lambda_{j}^{(1)}-{i\over 2}\over
\lambda+ \lambda_{j}^{(1)} +{i\over 2}} \;
{\lambda-\lambda_{j}^{(1)}-{i\over 2}\over \lambda-
\lambda_{j}^{(1)} +{i\over 2}}\,,
\non\\
A_{l}(\lambda) &=& \prod_{j=1}^{M^{(l)}}
{\lambda+\lambda_{j}^{(l)}+i+{il\over 2} \over \lambda+ \lambda_{j}^{(l)}
+{il\over2}} \; {\lambda-\lambda_{j}^{(l)}+{il\over 2}+i\over \lambda-
\lambda_{j}^{(l)} +{il\over 2}} \non\\
&& \times \prod_{j=1}^{M^{(l+1)}}{\lambda+ \lambda_{j}^{(l+1)}+{il\over
2}-{i\over 2}\over \lambda+ \lambda_{j}^{(l+1)} +{il\over 2} +{i\over 2}}\
{\lambda-\lambda_{j}^{(l+1)}+{il \over 2}-{i\over 2} \over
\lambda-\lambda_{j}^{(l+1)} + {il\over 2}+{i\over 2}} \,, \qquad l =
1,\ldots , k-1 \non\\
A_{k}(\lambda) &=& \prod_{j=1}^{M^{(k)}}
{\lambda+\lambda_{j}^{(k)}+{ik\over 2}-i\over \lambda+ \lambda_{j}^{(k)}
+{ik\over2}} \; {\lambda-\lambda_{j}^{(k)}+{ik\over 2}-i\over \lambda-
\lambda_{j}^{(k)} +{ik\over 2}} {\lambda+\lambda_{j}^{(k)}+{ik\over
2}+{i\over 2} \over \lambda+ \lambda_{j}^{(k)} +{ik\over 2}-{i\over 2}} \;
{\lambda-\lambda_{j}^{(k)}+ {ik\over 2}+{i\over 2}\over \lambda-
\lambda_{j}^{(k)} +{ik\over 2}-{i\over 2}}\,
\label{a3}
\end{eqnarray}
together with $A_{l}(\lambda) = A_{n-l-1}(-\lambda -i\kappa)$ for $l >k$.
We recall
that for $so(2k+1)$, $\kappa = k - \half$.

\subsection*{B. $\bf so(2k)$}
The dressing functions $A_l(\lambda)$
are the same as in the $so(2k+1)$ case for 
$ l=0,\ldots, k-3$, while
\begin{eqnarray}
A_{k-2}(\lambda) &=& \prod_{j=1}^{M^{(k-2)}}
{\lambda+\lambda_{j}^{(k-2)}+{ik\over 2} \over \lambda+ \lambda_{j}^{(k-2)}
+{ik\over2}-i} \; {\lambda-\lambda_{j}^{(k-2)}+{ik\over 2} \over \lambda-
\lambda_{j}^{(k-2)} +{ik\over 2}-i} \non\\
&& \times \prod_{j=1}^{M^{(+)}}{\lambda+\lambda_{j}^{(+)}+{ik\over
2}-{3i\over 2}\over \lambda+ \lambda_{j}^{(+)} +{ik\over 2}-{i\over 2}} \;
{\lambda-\lambda_{j}^{(+)}+{ik\over 2}-{3i\over 2}\over \lambda-
\lambda_{j}^{(+)} +{ik\over 2}-{i\over 2}} \non\\
&& \times \prod_{j=1}^{M^{(-)}}{\lambda+\lambda_{j}^{(-)}+{ik\over
2}-{3i\over 2}\over \lambda+ \lambda_{j}^{(-)} +{ik\over 2}-{i\over 2}} \;
{\lambda-\lambda_{j}^{(-)}+{ik\over 2}-{3i\over 2}\over \lambda-
\lambda_{j}^{(-)} +{ik\over 2}-{i\over 2}}\,, ~~ \non
\end{eqnarray}
\begin{eqnarray}
A_{k-1}(\lambda) &=& \prod_{j=1}^{M^{(+)}}
{\lambda+\lambda_{j}^{(+)}+{ik\over 2}-{3i\over 2}\over \lambda+
\lambda_{j}^{(+)} +{ik\over2}-{i \over 2}} \;
{\lambda-\lambda_{j}^{(+)}+{ik\over 2} -{3i \over 2}\over \lambda-
\lambda_{j}^{(+)} +{ik\over 2}-{i \over 2}} \non\\
&& \times \prod_{j=1}^{M^{(-)}} {\lambda+\lambda_{j}^{(-)}+{ik\over
2}+{i\over 2}\over \lambda+ \lambda_{j}^{(-)} + {ik\over 2}-{i\over 2}}
\; {\lambda-\lambda_{j}^{(-)}+{ik\over 2}+{i\over 2}\over \lambda-
\lambda_{j}^{(-)} +{ik\over 2}-{i\over 2}}\,
\label{a4}
\end{eqnarray}
with still $A_{l}(\lambda) = A_{n-l-1}(-\lambda -i\kappa)$ for $l >k-1$, with
$\kappa = k-1$.

\subsection*{C. $\bf sp(2k)$}
Similarly, the dressing functions are the same as in
the $so(2k+1)$ case for $ l=0, \ldots, k-3$, and:
\begin{eqnarray}
A_{k-2}(\lambda) &=& \prod_{j=1}^{M^{(k-2)}}
{\lambda+\lambda_{j}^{(k-2)}+{ik\over 2} \over
\lambda+\lambda_{j}^{(k-2)}+{ik\over2}-i} \;
{\lambda-\lambda_{j}^{(k-2)}+{ik\over 2} \over \lambda- \lambda_{j}^{(k-2)}
+{ik\over 2}-i} \non\\
&& \times \prod_{j=1}^{M^{(k-1)}}{\lambda+\lambda_{j}^{(k-1)}+{ik\over
2}-{3i\over 2}\over \lambda+ \lambda_{j}^{(k-1)} +{ik\over 2}-{i\over2}} \;
{\lambda-\lambda_{j}^{(k-1)}+{ik\over 2}-{3i\over 2}\over \lambda-
\lambda_{j}^{(k-1)} +{ik\over 2}-{i\over 2}} \non\\
A_{k-1}(\lambda) &=&
\prod_{j=1}^{M^{(k-1)}}{\lambda+\lambda_{j}^{(k-1)}+{ik\over 2}+ {i\over
2}\over \lambda+ \lambda_{j}^{(k-1)} +{ik\over2}-{i \over 2}} \; {\lambda-
\lambda_{j}^{(k-1)}+{ik\over 2}+{i\over 2}\over \lambda-
\lambda_{j}^{(k-1)} + {ik\over 2}-{i \over 2}} \non\\
&& \times \prod_{j=1}^{M^{(k)}}{\lambda+\lambda_{j}^{(k)}+{ik\over 2}
-{3i\over 2}\over \lambda+ \lambda_{j}^{(k)} +{ik\over 2}+{i\over 2}} \;
{\lambda- \lambda_{j}^{(k)}+{ik\over 2}-{3i\over 2}\over \lambda-
\lambda_{j}^{(k)} + {ik\over 2}+{i\over 2}}\, ~~\label{a5}
\end{eqnarray}
In addition $A_{l}(\lambda) = A_{n-l-1}(-\lambda -i\kappa)$ for $l >k-1$;
note that now $\kappa = k+1$.

\subsection{Bethe Ansatz equations for $K^-=1$}
The above dressing functions $A_{l}$ satisfy all the imposed constraints
and they are unambiguously defined. \\
Requiring now the analyticity of the eigenvalues, we deduce the Bethe
Ansatz equations. More specifically, successive $A_{l}$'s have common
poles, which must disappear. Hence, the sum of corresponding residues of
$A_l$ and $A_{l+1}$ in the eigenvalue expression (\ref{eigen}) must be
zero. The Bethe Ansatz equations immediately follow from this condition.

Let us define
\begin{eqnarray}
e_{x}(\lambda) = {\lambda +{ix \over 2} \over \lambda -{ix \over 2}},
\end{eqnarray}
for any $x$. Then, the Bethe Ansatz equations read:

\subsection*{A. $\bf so(2k+1)$}
\begin{eqnarray}
e_{1}(\lambda_{i}^{(1)})^{2N} &\!\!=\!\!& \prod_{j=1,j \ne i}^{M^{(1)}}
e_{2}(\lambda_{i}^{(1)} - \lambda_{j}^{(1)})\ e_{2}(\lambda_{i}^{(1)} +
\lambda_{j}^{(1)})\ \prod_{ j=1}^{M^{(2)}}e_{-1}(\lambda_{i}^{(1)} -
\lambda_{j}^{(2)})\ e_{-1}(\lambda_{i}^{(1)} + \lambda_{j}^{(2)})\,, \non\\
1 &\!\!=\!\!& \prod_{j=1,j \ne i}^{M^{(l)}} e_{2}(\lambda_{i}^{(l)} -
\lambda_{j}^{(l)})\ e_{2}(\lambda_{i}^{(l)} + \lambda_{j}^{(l)})\ \prod_{
\eta = \pm 1}\prod_{ j=1}^{M^{(l+\eta)}}e_{-1}(\lambda_{i}^{(l)} -
\lambda_{j}^{(l+\eta)})\ e_{-1}(\lambda_{i}^{(l)} + \lambda_{j}^{(l+\eta)})
\non\\
&& l= 2, \ldots , k-1 \non\\
1 &\!\!=\!\!& \prod_{j=1,j \ne i}^{M^{(k)}} e_{1}(\lambda_{i}^{(k)} -
\lambda_{j}^{(k)})\ e_{1}(\lambda_{i}^{(k)} + \lambda_{j}^{(k)})\ \prod_{
j=1}^{M^{(k-1)}}e_{-1}(\lambda_{i}^{(k)} - \lambda_{j}^{(k-1)})\
e_{-1}(\lambda_{i}^{(k)} + \lambda_{j}^{(k-1)}) \non\\
\label{BAE1}
\end{eqnarray}

\subsection*{B. $\bf so(2k)$}
The first $k-3$ equations are the same as in $so(2k+1)$, see eq.
(\ref{BAE1}), but the last three equations are modified, namely
\begin{eqnarray}
1 &\!\!=\!\!& \prod_{j=1,j \ne i}^{M^{(k-2)}} e_{2}(\lambda_{i}^{(k-2)} -
\lambda_{j}^{(k-2)})\ e_{2}(\lambda_{i}^{(k-2)} + \lambda_{j}^{(k-2)})\
\prod_{ j=1}^{M^{(k-3)}}e_{-1}(\lambda_{i}^{(k-2)} - \lambda_{j}^{(k-3)})\
e_{-1}(\lambda_{i}^{(k-2)} + \lambda_{j}^{(k-3)})\non\\
&& \times \prod_{ j=1}^{M^{(+)}}e_{-1}(\lambda_{i}^{(k-2)} -
\lambda_{j}^{(+)})\ e_{-1}(\lambda_{i}^{(k-2)} + \lambda_{j}^{(+)})
\ \prod_{ j=1}^{M^{(-)}}e_{-1}(\lambda_{i}^{(k-2)} -
\lambda_{j}^{(-)})\ e_{-1}(\lambda_{i}^{(k-2)} + \lambda_{j}^{(-)})
\non\\
1 &\!\!=\!\!& \prod_{j=1,j \ne i}^{M^{(\tau)}} e_{2}(\lambda_{i}^{(\tau)} -
\lambda_{j}^{(\tau)})\ e_{2}(\lambda_{i}^{(\tau)} + \lambda_{j}^{(\tau)})\
\prod_{ j=1}^{M^{(k-2)}}e_{-1}(\lambda_{i}^{(\tau)} - \lambda_{j}^{(k-2)})\
e_{-1}(\lambda_{i}^{(\tau)} + \lambda_{j}^{(k-2)}),\ \tau=\pm.
\label{BAE2}
\end{eqnarray}

\subsection*{C. $\bf sp(2k)$}
The first $k-2$ equations are the same as in the $so(2k+1)$, while the last
two equations are given by
\begin{eqnarray}
1 &\!\!=\!\!& \prod_{j=1,j \ne i}^{M^{(k-1)}} e_{2}(\lambda_{i}^{(k-1)}
- \lambda_{j}^{(k-1)})\ e_{2}(\lambda_{i}^{(k-1)} +
\lambda_{j}^{(k-1)})\  \prod_{
j=1}^{M^{(k-2)}}e_{-1}(\lambda_{i}^{(k-1)} - \lambda_{j}^{(k-2)})\
e_{-1}(\lambda_{i}^{(k-1)} + \lambda_{j}^{(k-2)})\non\\
&& \times \prod_{ j=1}^{M^{(k)}}e_{-2}(\lambda_{i}^{(k-1)} -
\lambda_{j}^{(k)})\ e_{-2}(\lambda_{i}^{(k-1)} + \lambda_{j}^{(k)}) \non\\
1 &\!\!=\!\!& \prod_{j=1,j \ne i}^{M^{(k)}} e_{4}(\lambda_{i}^{(k)} -
\lambda_{j}^{(k)})\ e_{4}(\lambda_{i}^{(k)} + \lambda_{j}^{(k)})\ \prod_{
j=1}^{M^{(k-1)}}e_{-2}(\lambda_{i}^{(k)} - \lambda_{j}^{(k-1)})\
e_{-2}(\lambda_{i}^{(k)} + \lambda_{j}^{(k-1)}).
\label{BAE3}
\end{eqnarray}

\bigskip
For all these cases, we recover the rational limits of the equations given in
\cite{AMN1,AMN2}. 

\subsection{Eigenvalues and Bethe Ansatz equations for diagonal $K^-$}

Until now, we have considered the case of trivial boundary effects
$K^+=K^-=1$. We here come to the main point of our derivation and insert
non-trivial boundary effects. We shall then rederive modified Bethe Ansatz
equations. We choose $K^{-}$ to be one of the diagonal solutions D1, D2,
D3, D4. We consider, for simplicity but without loss of generality, $K^{+}
=1$. Remark that the pseudo-vacuum remains an exact eigenstate after
this modification.

Let us rewrite the solutions D1, D2, D3 and D4 in a slightly modified
notation, which we are going to use from now on.

\medskip

{\bf D1:} The solution D1 can be written in the following form
\begin{eqnarray}
K(\lambda) = diag( \alpha, \ldots ,\alpha, \beta, \dots, \beta) \,.
\label{1}
\end{eqnarray}
The number of $\alpha 's$ is equal to the number of $\beta 's$, so that
this solution exists only for the $so(2k)$ and $sp(2k)$ cases as stated in
Proposition \ref{prop:D}, and
\begin{eqnarray}
\alpha(\lambda) = -\lambda +i\xi, ~~\beta(\lambda) = \lambda +i\xi,
\end{eqnarray}
where $\xi ={1\over c}$ is the free
boundary parameter.

\medskip

{\bf D2:} Solution D2 can be written as
\begin{eqnarray}
K(\lambda) = diag( \alpha, \beta, \ldots, \beta,  \gamma) \label{2}
\end{eqnarray}
with
\begin{eqnarray}
\alpha(\lambda) = { -\lambda +i\xi_{1} \over \lambda +i\xi_{1}},
~~\beta(\lambda) = 1, ~~\gamma(\lambda)
= { -\lambda +i\xi_{n} \over \lambda +i\xi_{n}} \,,
\end{eqnarray}
where $\xi_{1} = - {1\over c_{1}}$, $\xi_{n}= -{1\over c_{n}} $ are the
boundary parameters which satisfy the constraint
\begin{eqnarray}
\xi_{1} +\xi_{n} = \kappa -1.
\end{eqnarray}
We remind that this solution exists for $so(n)$ algebras, but not for
$sp(n)$.

\medskip

{\bf D3:} Solution D3 has the form
\begin{eqnarray}
K(\lambda) = diag( \alpha, \ldots ,\alpha, \beta, \dots, \beta, \alpha,
\ldots ,\alpha).
\label{3}
\end{eqnarray}
The number of $\alpha 's$ is $2m$ and the number of $\beta 's$ is $n-2m$,
where the free integer parameter $m$ obeys 
$1 \le m \le k-1$ for $so(2k)$, $sp(2k)$ and $1 \le m \le k$ for
$so(2k+1)$. Moreover
\begin{eqnarray}
\alpha(\lambda) = -\lambda +i\xi, ~~\beta(\lambda) = \lambda +i\xi
\end{eqnarray}
where $\xi ={n\over 4} -m$ has a fixed value, unlike the $su(n)$ case where
the corresponding solution has a free boundary parameter (see
\cite{dVG,done}).

\medskip

{\bf D4:} Finally solution D4, which holds only for the $so(4)$ case, can
be written in the following form
\begin{eqnarray}
K(\lambda) = diag( \alpha,\beta,  \gamma, \delta) \label{4}
\end{eqnarray}
where
\begin{eqnarray}
\alpha(\lambda) &=& (-\lambda +i\xi_{-})(-\lambda +i\xi_{+}),
~~\beta(\lambda) = (\lambda +i\xi_{-})
(-\lambda +i\xi_{+}) \,, \non\\
\gamma(\lambda) &=& (-\lambda +i\xi_{-})(\lambda +i\xi_{+}),
~~\delta(\lambda) = (\lambda +i\xi_{-})(\lambda +i\xi_{+}) \,.
\end{eqnarray}
$\xi_{-} = {1\over c_{2}}$, $\xi_{+} = {1\over c_{3}}$ are both free
parameters.

\null

We now come to the explicit expression of the eigenvalues when $K^-$ is one
of the above mentioned solutions. We should point out that the dressing
functions are related to the bulk behavior of the chain and thus they are
form-invariant under changes of boundary conditions. Indeed what is
modified in the expression of the eigenvalues (\ref{eigen}) are the $g_{l}$
functions which characterise the boundary effects. We call the new $g_{l}$
functions $\tilde g_{l}$.


\medskip

{\bf D1:} As already mentioned, the solution D1 can only be applied in
$so(2k)$ and $sp(2k)$. In this case we have
\begin{eqnarray}
\tilde g_{l}(\lambda)&=&(-\lambda + i\xi)g_{l}(\lambda), ~~l=0,\ldots ,k-1
\non\\
\tilde g_{l}(\lambda)&=&(\lambda + i\xi +i\kappa)g_{l}(\lambda),
~~l=k,\ldots ,2k-1 \label{tg1}
\end{eqnarray}
where $g(\lambda)$ are given by (\ref{g1})--(\ref{g3}). The system with
such boundaries has a residual symmetry $sl(k)$ in both cases, which
immediately follows from the structure of the corresponding $K$ matrix.

\medskip

{\bf D2:} We have
\begin{eqnarray}
\tilde g_{0}(\lambda)&=&{(-\lambda + i\xi_{1}) \over (\lambda + i\xi_{1})}
g_{0}(\lambda),
~~\tilde g_{n-1}(\lambda)={(\lambda + i\xi_{1} +i)\over (\lambda +
i\xi_{1})}{(\lambda + i\xi_{1}+i\kappa)
\over (-\lambda -i\kappa+ i\xi_{1}+i)}g_{n-1}(\lambda), \non\\  \tilde
g_{l}(\lambda)&=&{(\lambda + i\xi_{1}+i)
\over (\lambda + i\xi_{1})} g_{l}(\lambda), ~~l=1,\ldots ,n-2 \label{tg2}
\end{eqnarray}
Again the $g_{l}(\lambda)$ are given by (\ref{g1})--(\ref{g3}). Similarly,
from the structure of the $K$ matrix we conclude that the residual symmetry
is $so(n-2) \otimes so(2)$.

\medskip

{\bf D3:} For the D3 solution we find the following modified $g$ functions
\begin{eqnarray}
\tilde g_{l}(\lambda)&=& (-\lambda +i\xi)g_{l}(\lambda), ~~l=0,\ldots ,m-1
\non\\
\tilde g_{l}(\lambda) &=&(\lambda + {i\kappa \over 2} \pm {i\over
2})g_{l}(\lambda), ~~l=m,\ldots ,n-m-1 \non\\
\tilde g_{l}(\lambda) &=& (-\lambda -i\kappa - i\xi) {(\lambda + {i\kappa
\over 2} \pm{i\over 2}) \over (\lambda + {i\kappa \over 2} \mp{i\over 2})}
g_{l}(\lambda), ~~l=n-m,\ldots ,n-1. \label{tg3}
\end{eqnarray}
The symmetry of the transfer matrix for this $K$ matrix is $so(n-2m)
\otimes so(2m)$, (resp. $sp(n-2m) \otimes sp(2m)$).

\medskip

{\bf D4:} Finally for solution D4 the modified $g$ functions are given by
\begin{equation}
\begin{array}{llllll}
\displaystyle\tilde g_{0}(\lambda)&=& (-\lambda +i\xi_{-})(-\lambda
+i\xi_{+})g_{0}(\lambda),\qquad & \tilde g_{1}(\lambda)&=&(\lambda
+i\xi_{-}+i)(-\lambda +i\xi_{+})g_{1}(\lambda) \non\\
\displaystyle\tilde g_{2}(\lambda)&=& (\lambda +i\xi_{+}+i)(-\lambda
+i\xi_{-})g_{2}(\lambda), & \tilde g_{3}(\lambda)&=&(\lambda
+i\xi_{-}+i)(\lambda +i\xi_{+}+i)g_{3}(\lambda).
\end{array}
\label{tg4}
\end{equation}
We now formulate the Bethe Ansatz equations for the general diagonal
solutions. The only modifications induced on equations
(\ref{BAE1})--(\ref{BAE3}) are the following for each solution:
\begin{itemize}
\item[{\bf D1}]
The factor $-e_{2\xi+\kappa}^{-1}(\lambda)$ appears in the LHS of the
$k^{th}$ Bethe equation.
\item[{\bf D2}]
The factor $-e_{2\xi_{1}+1}^{-1}(\lambda)$ appears in the LHS of the first
Bethe equation.
\item[{\bf D3}]
The factor $-e_{2\xi+m}^{-1}(\lambda)$ appears in the LHS of the $m^{th}$
Bethe equation with $m=1,\ldots ,k-1$ for $so(2k+1)$, $sp(2k)$
and $m=1,\ldots ,k-2$ for $so(2k)$. 
For $so(2k)$,  when $m=k-1$, 
the factor $-e_{1}^{-1}(\lambda)$ appears in the LHS
of the $(k-1)^{th}$ and $k^{th}$ Bethe Ansatz equations.
\end{itemize}
We treat solution D4 separately in the next section.

\section{ Ground state and excitations}
\setcounter{equation}{0}

The next step is to determine the ground state and the low-lying
excitations of the model. One of the main aims of this work is indeed the
computation of the scattering of the low-lying excitations off the
boundaries. The bulk scattering for these models has been already studied
in \cite{ogiev}, nevertheless we are going to rederive these results as a
check in the next section.

We recall that the quantum numbers that describe a state are given by
(\ref{qn}), and the energy is derived via the relation $H={d\over d
\lambda}t(\lambda) \vert_ {\lambda = 0}$. It is given by
\begin{eqnarray}
E = -{1\over 2 \pi } \sum_{j=1}^{M^{(1)}} {1\over
{\left(\lambda_{j}^{(1)}\right)}^{2} + {1\over 4}} \,.
\label{energy}
\end{eqnarray}
In what follows we write the Bethe Ansatz equations for the ground state
and the low-lying excitations (holes) of the models under study. Bethe
Ansatz equations may in general only be solved in the thermodynamic limit
$N\to \infty$. In this limit, it is assumed that a state is described in
particular by the \emph{density functions} $\sigma^{l}(\lambda)$ of the
parameters $\lambda_{i}^{(l)}$. \\
We here make the hypothesis that non trivial boundary effects do not
modify the nature of the ground state and excited states but only the
values of the Bethe parameters.

\subsection*{A. $ \bf so(2k+1)$}
Let us first consider the $so(2k+1)$ case, for which the ground state
consists of $k-1$ filled Dirac seas, whereas the $k^{th}$ sea is filled
with two-strings of the form $\lambda_{0} \pm {i\over 4}$. Remark the shift
${i \over 4}$, instead of ${i \over 2}$ usually expected, because the
length of the $k^{th}$ root in the $so(2k+1)$ case is half the length of
the other roots (see also \cite{ogiev}). This leads us to rewrite the Bethe
Ansatz equations for the ground state and the low-lying excitations (holes)
in the following form, inserting the two-strings contribution in the
$k^{th}$ set:
\begin{eqnarray}
e_{1}(\lambda_{i}^{(1)})^{2N+1} &\!\!=\!\!& -\prod_{j=1}^{M^{(1)}}
e_{2}(\lambda_{i}^{(1)} - \lambda_{j}^{(1)})\ e_{2}(\lambda_{i}^{(1)} +
\lambda_{j}^{(1)})\ \prod_{j=1}^{M^{(2)}}e_{-1}(\lambda_{i}^{(1)} -
\lambda_{j}^{(2)})\ e_{-1}(\lambda_{i}^{(1)} + \lambda_{j}^{(2)})\,, \non\\
e_{1}(\lambda_{i}^{(l)}) &\!\!=\!\!& -\prod_{j=1}^{M^{(l)}}
e_{2}(\lambda_{i}^{(l)} - \lambda_{j}^{(l)})\ e_{2}(\lambda_{i}^{(l)} +
\lambda_{j}^{(l)})\ \prod_{\tau =\pm 1}
\prod_{j=1}^{M^{(l+\tau)}}e_{-1}(\lambda_{i}^{(l)} -
\lambda_{j}^{(l+\tau)})\ e_{-1}(\lambda_{i}^{(l)} + \lambda_{j}^{(l+\tau)})
\non\\
&& (l= 2, \ldots , k-2) \non\\
e_{1}(\lambda_{i}^{(k-1)}) &\!\!=\!\!&
-\prod_{j=1}^{M^{(k-1)}}e_{2}(\lambda_{i}^{(k-1)} - \lambda_{j}^{(k-1)})\
e_{2}(\lambda_{i}^{(k-1)} + \lambda_{j}^{(k-1)}) \non\\
&& \times \prod_{j=1}^{M^{(k-2)}}e_{-1}(\lambda_{i}^{(k-1)} -
\lambda_{j}^{(k-2)})\ e_{-1}(\lambda_{i}^{(k-1)} + \lambda_{j}^{(k-2)})\
\non\\
&& \times \prod_{j=1}^{M^{(k)}}e_{-{1\over 2}}(\lambda_{i}^{(k-1)} -
\lambda_{j}^{(k)}) \ e_{-{1\over 2}}(\lambda_{i}^{(k-1)} +
\lambda_{j}^{(k)}) e_{-{3\over 2}}(\lambda_{i}^{(k-1)} -
\lambda_{j}^{(k)})\ e_{-{3\over 2}}(\lambda_{i}^{(k-1)} +
\lambda_{j}^{(k)}) \non\\
e_{1}(\lambda_{i}^{(k)}) &\!\!=\!\!& -\prod_{j=1}^{M^{(k)}}
e_{1}(\lambda_{i}^{(k)} - \lambda_{j}^{(k)})^{2} \
e_{1}(\lambda_{i}^{(k)} + \lambda_{j}^{(k)})^{2} \ e_{2}(\lambda_{i}^{(k)} -
\lambda_{j}^{(k)})\ e_{2}(\lambda_{i}^{(k)} + \lambda_{j}^{(k)}) \non \\
&& \times \prod_{j=1}^{M^{(k-1)}}e_{-{1\over 2}}(\lambda_{i}^{(k)} -
\lambda_{j}^{(k-1)})\ e_{-{1\over 2}}(\lambda_{i}^{(k)} +
\lambda_{j}^{(k-1)}) e_{-{3\over 2}}(\lambda_{i}^{(k)} -
\lambda_{j}^{(k-1)})\ e_{-{3\over 2}}(\lambda_{i}^{(k)} +
\lambda_{j}^{(k-1)}) \non\\
\label{BAE1g}
\end{eqnarray}
For the general diagonal solutions of the reflection equation we have to
multiply: 
\begin{itemize}
\item[--] {\bf for D2}, the LHS of the $1^{st}$ Bethe Ansatz equation with
$-e_{2\xi_{1}
+1}^{-1}$; 
\item[--] {\bf for D3}, we multiply the LHS of the $m^{th}$ Bethe Ansatz
equation with
$-e_{2\xi +m}^{-1}$, $1 \le m \le k-1$. 
\end{itemize}
We are interested in the low-lying excitations, which are holes in the
filled sea and are highest weight representations of $so(2k+1)$. We
restrict ourselves here to the states with $\nu^{(l)}$ holes in the $l$
sea, which correspond to the vector representations of $so(2k+1)$, and also
to holes in the $k$ sea, which corresponds to the $2^{k}$-dimensional
spinor representation (see also \cite{reshe}). We convert the sums into
integrals by employing the following approximate relation \cite{GMN}
\begin{eqnarray}
{1\over N} \sum_{i=1}^{M} f(\lambda_{i}^{(l)})= \int_{0}^{\infty} d\lambda
f(\lambda) \sigma^{l}(\lambda) - {1\over N} \sum_{i=1}^{\nu^{(l)}} f(\tilde
\lambda_{i}^{(l)}) -{1\over 2N}f(0)+O(\frac{1}{N^2})
\label{mac}
\end{eqnarray}
where the correction terms take into account the $\nu^{(l)}$ holes located
at values $\tilde \lambda_{i}^{(l)}$ and the contribution at $0^+$. We
shall denote by $\hat{f}(\omega)$ the Fourier transform of any function
$f(\lambda)$. \\
Once we take the logarithm and the derivative of (\ref{BAE1g}), we extract
the densities from the equation
\begin{eqnarray}
\hat {\cal K}(\omega) \hat \sigma(\omega) = \hat a(\omega) +{1\over N} \hat
F(\omega) +{1\over N} \hat G(\omega, \xi) \label{sigma}
\end{eqnarray}
where $\displaystyle a_{x}(\lambda) = \frac{i}{2\pi} \; \frac{d}{d\lambda}
\, \ln e_{x}(\lambda)$ and $\hat a_{x}(\omega) =e^{-{x \omega \over 2}}$.
We have introduced
\begin{eqnarray}
a(\lambda) = \left
(\begin{array}{c}
2a_{1}(\lambda) \\
0 \\
\vdots \\
0 \\
\end{array}
\right)\,, ~~\sigma(\lambda) = \left (
\begin{array}{c}
\sigma^{1}(\lambda) \\
\vdots  \\
\sigma^{l}(\lambda)\\
\vdots \\
\sigma^{k}(\lambda)  \\
\end{array}
\right)\,.
\label{cols}
\end{eqnarray}
$F(\lambda)$, $G(\lambda, \xi)$ are $k$ component vectors as well with
\begin{eqnarray}
F^{j}(\lambda) &\!\!=\!\!& a_{1}(\lambda)\delta_{j1}-a_{1}(\lambda)
+a_{2}(\lambda)+ \sum_{j=1}^{\nu^{(l)}} \Big(a_{2}(\lambda -\tilde
\lambda_{j}^{(l)})+ a_{2}(\lambda +\tilde \lambda_{j}^{(l)})\Big)
\delta_{lj} \non\\
&-&\sum_{j=1}^{\nu^{(l)}} \Big(a_{1}(\lambda -\tilde
\lambda_{j}^{(l)})+a_{1}(\lambda +\tilde \lambda_{j}^{(l)})\Big)
(\delta_{j,l+1} +\delta_{j,l-1}) \,, \qquad (j=1,\ldots,k-2) \non\\
F^{k-1}(\lambda) &\!\!=\!\!& a_{2}(\lambda)-(a_{{1\over 2}}(\lambda)+a_{{3
\over 2}}(\lambda))- \sum_{j=1}^{\nu^{(l)}} \Big(a_{1}(\lambda -\tilde
\lambda_{j}^{(l)})+a_{1}(\lambda +\tilde \lambda_{j}^{(l)})\Big)
\delta_{l,k-2} \non\\
&+& \sum_{j=1}^{\nu^{(l)}} \Big(a_{2}(\lambda -\tilde
\lambda_{j}^{(l)})+a_{2}(\lambda +\tilde \lambda_{j}^{(l)})\Big)
\delta_{l,k-1} -\sum_{j=1}^{\nu^{(l)}} \Big((a_{{1\over 2}}+a_{{3 \over
2}})(\lambda-\tilde \lambda_{j}^{(l)})+(a_{{1\over 2}}+a_{{3 \over
2}})(\lambda+\tilde \lambda_{j}^{(l)})\Big) \delta_{kl} \non\\
F^{k}(\lambda) &\!\!=\!\!& 3a_{1}(\lambda) + a_{2}(\lambda)-(a_{{1\over
2}}(\lambda)+a_{{3 \over 2}}(\lambda))-\sum_{j=1}^{\nu^{(l)}}
\Big((a_{{1\over 2}}+a_{{3 \over 2}})(\lambda-\tilde
\lambda_{j}^{(l)})+(a_{{1\over 2}}+a_{{3 \over 2}})(\lambda+\tilde
\lambda_{j}^{(l)})\Big) \delta_{k-1,l} \non\\
&+&\sum_{j=1}^{\nu^{(l)}} \Big((2a_{1}+a_{2})(\lambda-\tilde
\lambda_{j}^{(l)})+(2a_{1}+a_{2})(\lambda+\tilde \lambda_{j}^{(l)})\Big)
\delta_{kl}
\label{f1}
\end{eqnarray}
and the $\xi$ dependent part is
\begin{eqnarray}
\mbox{\bf for D2:}\qquad ~~~ G^{j}(\lambda, \xi) &\!\!=\!\!&
-a_{2\xi_{1}+1}(\lambda)\delta_{j1}, \non\\
\mbox{\bf for D3:}\qquad ~~~ G^{j}(\lambda, \xi) &\!\!=\!\!&
-a_{2\xi+m}(\lambda)\delta_{jm} ~~(m=1,...,k-1).
\label{xi1}
\end{eqnarray}
Finally, the non vanishing entries of $\hat{\cal K}$ are
\begin{eqnarray}
\hat {\cal K}_{ij}(\omega) &\!\!=\!\!& (1+\hat a_{2}(\omega))\delta_{ij} -
\hat a_{1}(\omega)(\delta_{i,j+1}+\delta_{i,j-1}), ~~i,j =1,\ldots, k-1
\non\\
\hat {\cal K}_{k-1,k}(\omega) &\!\!=\!\!& \hat {\cal K}_{k,k-1}(\omega) =
-(\hat a_{{1\over 2}}(\omega) +\hat a_{{3\over 2}}(\omega)), \non\\
\hat {\cal K}_{kk}(\omega) &\!\!=\!\!& 1+2\hat a_{1}(\omega)+\hat
a_{2}(\omega).
\label{K1}
\end{eqnarray}
The solution of (\ref{sigma}) has the form
\begin{eqnarray}
\sigma(\lambda) =
2 \epsilon(\lambda) +{1\over N} \Phi_{0}(\lambda)+{1\over N}
\Phi_{1}(\lambda, \xi) \label{sigma2}
\end{eqnarray}
where $\epsilon$ and $\Phi_{0,1}$ are $k$ component vectors with
\begin{eqnarray}
\hat  \epsilon^{i}(\omega) = \hat {\cal R}_{i1}(\omega) \hat
a_{1}(\omega), ~~\hat \Phi_{0}^{i}(\omega) =\sum_{j=1}^{k}\hat
{\cal R}_{ij}(\omega) \hat F^{j}(\omega), ~~\hat \Phi_{1}^{i}(\omega,\xi)
=\sum_{j=1}^{k}\hat {\cal R}_{ij}(\omega) \hat G^{j}(\omega,\xi)
\label{ef1}
\end{eqnarray}
$\hat{\cal R} = \hat{\cal K}^{-1}$ and $\hat\epsilon^{j}$ is the energy of
a hole in 
the $j$ sea which can be written in terms of hyperbolic functions
\begin{eqnarray}
\hat \epsilon^{j}(\omega) ={\cosh(k-{1\over 2} -j){\omega \over 2} \over
\cosh(k-{1\over 2}){\omega \over 2}}, ~~j=1,\ldots,k-1, ~~~\hat
\epsilon^{k}(\omega) ={1 \over 2\cosh(k-{1\over 2}){\omega \over 2}}
\label{ener1}
\end{eqnarray}
\begin{eqnarray}
&& \hat {\cal R}_{ij}(\omega) =e^{{\omega \over 2}} {\sinh \min(i,j) {\omega
\over 2} \cosh(k-{1\over 2} -\max(i,j)){\omega \over 2} \over \cosh
(k-{1\over 2}){\omega \over 2}\sinh{\omega \over 2} }, ~~i,j = 1, \ldots
,k-1 \non\\
&& \hat {\cal R}_{jk}(\omega)=\hat {\cal R}_{kj}(\omega) ={e^{{\omega \over
2}}\over 2}{\sinh {j\omega \over 2} \over \cosh (k-{1\over 2}){\omega \over
2}\sinh{\omega \over 2}}, ~~j = 1, \ldots ,k-1 \non\\
&& \hat {\cal R}_{kk}(\omega)={e^{{\omega \over 2}}\over 2}{\sinh {k\omega
\over 2} \over 2\cosh{\omega \over 4} \cosh (k-{1\over 2}){\omega \over
2}\sinh{\omega \over 2}}
\label{RR1}
\end{eqnarray}

\subsection*{B. $ \bf so(2k)$}
In this case, the ground state consists of $k$ filled Dirac seas of real
strings. Therefore the Bethe Ansatz equations have exactly the same form as
in (\ref{BAE2}). For the general diagonal solutions of the reflection
equation, we have to multiply:
\begin{itemize}
\item[--]
{\bf for D1}, the LHS of the $k^{th}$ Bethe Ansatz equation with $-e_{2\xi
+\kappa}^{-1}$,
\item[--]
{\bf for D2}, the LHS of the $1^{st}$ equation with $-e_{2\xi_{1}+1}^{-1}$,
\item[--]
 {\bf for D3}, the LHS of the $m^{th}$ equation with
 $-e_{2\xi+m}^{-1}$ for $m=1,\ldots,k-2$, while the LHS of the
 $(k-1)^{th}$ and $k^{th}$ equations for $m=k-1$ are 
 multiplied by $-e_{1}^{-1}$.
\end{itemize}
In this case as well, we restrict ourselves to states with $\nu^{(l)}$
holes in the $l$ sea. Note that now the spinor representation splits into
two spinor representations of dimension $2^{k-1}$ (see also \cite{reshe}),
and the holes in the $+$, $-$ sea correspond exactly to these two spinor
representations. The densities satisfy the same equation (\ref{sigma}) as
in the $so(2k+1)$ case with $\sigma$ and $a$ given by (\ref{cols}) and
\begin{eqnarray}
\hat \epsilon^{j}(\omega)
={\cosh(k-1 -j){\omega \over 2} \over \cosh(k-1){\omega \over 2}},
~~j=1,\ldots,k-2 ~~\mbox{and}~~ \hat \epsilon^{\pm}(\omega) ={1 \over
2\cosh(k-1){\omega \over 2}}\label{ener2}
\end{eqnarray}
\begin{eqnarray}
F^{j}(\lambda) &\!\!=\!\!& a_{1}(\lambda)\delta_{j1}-a_{1}(\lambda)
+a_{2}(\lambda)+ \sum_{j=1}^{\nu^{(l)}} \Big(a_{2}(\lambda -\tilde
\lambda_{j}^{(l)})+a_{2}(\lambda +\tilde \lambda_{j}^{(l)})\Big)
\delta_{lj} \non\\
&-&\sum_{j=1}^{\nu^{(l)}} \Big(a_{1}(\lambda -\tilde
\lambda_{j}^{(l)})+a_{1}(\lambda +\tilde \lambda_{j}^{(l)})\Big)
(\delta_{j,l+1} +\delta_{j,l-1}), \qquad (j=1,\ldots ,k-3) \non\\
F^{k-2}(\lambda) &\!\!=\!\!& a_{2}(\lambda)-2a_{1}(\lambda) +
\sum_{j=1}^{\nu^{(l)}} \Big(a_{2}(\lambda -\tilde
\lambda_{j}^{(l)})+a_{2}(\lambda +\tilde \lambda_{j}^{(l)})\Big)
\delta_{l,k-2} \non\\
&-&\sum_{j=1}^{\nu^{(l)}} \Big(a_{1}(\lambda -\tilde
\lambda_{j}^{(l)})+a_{1}(\lambda +\tilde \lambda_{j}^{(l)})\Big)
(\delta_{l,k-3} +\delta_{l+}+\delta_{l-}) \non\\
F^{\pm}(\lambda) &\!\!=\!\!& a_{2}(\lambda) + \sum_{j=1}^{\nu^{(l)}}\Big(a_{2}
(\lambda -\tilde \lambda_{j}^{(l)})+a_{2}(\lambda +\tilde
\lambda_{j}^{(l)})\Big) \delta_{l\pm}-\sum_{j=1}^{\nu^{(l)}}
\Big(a_{1}(\lambda -\tilde \lambda_{j}^{(l)})+a_{1}(\lambda +\tilde
\lambda_{j}^{(l)})\Big) \delta_{l,k-2} \non \\
\label{f2}
\end{eqnarray}
\begin{eqnarray}
\mbox{-- \bf for D1:} && G^{j}(\lambda, \xi) =
-a_{2\xi+\kappa}(\lambda)\delta_{jk}, \non\\
\mbox{-- \bf for D2:} && G^{j}(\lambda, \xi) =
-a_{(2\xi_{1}+1)}(\lambda)\delta_{j1}, \non\\
\mbox{-- \bf for D3:} &&\left\{\begin{array}{l} G^{j}(\lambda, \xi) =
-a_{(2\xi+m)}(\lambda)\delta_{jm} ~~(m=1, \ldots , k-2), \\[.7ex]
G^{i}(\lambda)=
-a_{1}(\lambda)(\delta_{i,k-1}+\delta_{i,k}) ~~(m=k-1) \;.
\end{array}\right. 
\label{xi2}
\end{eqnarray}
Again, the non-vanishing entries of $\hat {\cal K}$ are
\begin{eqnarray}
&& \hat {\cal K}_{ij}(\omega) = (1+\hat a_{2}(\omega))\delta_{ij} - \hat
a_{1}(\omega)(\delta_{i,j+1}+\delta_{i,j-1}), ~~i,j =1,\ldots, k-2, \non\\
&& \hat {\cal K}_{k-2,\pm}(\omega) = \hat {\cal K}_{\pm, k-2}(\omega) =
-\hat a_{1}(\omega), \non\\
&& \hat {\cal K}_{--}(\omega)=\hat {\cal K}_{++}(\omega) = 1+\hat
a_{2}(\omega) \;\mbox{ and }\; \hat {\cal K}_{+-}(\omega)=\hat {\cal
K}_{-+}(\omega)=0.
\label{K2}
\end{eqnarray}
We solve equation (\ref{sigma}) and find the densities $\sigma^{i}$ which
describe a Bethe Ansatz state. The solution of (\ref{sigma}) has the same
form as in (\ref{sigma2}) with
\begin{eqnarray}
&&\hat {\cal R}_{ij}(\omega) =e^{{\omega \over 2}} {\sinh \min(i,j) {\omega
\over 2} \cosh(k-1 -\max(i,j)){\omega \over 2} \over \cosh (k-1){\omega
\over 2}\sinh{\omega \over 2} }, ~~i,j = 1, \ldots ,k-2, \non\\
&& \hat {\cal R}_{j\pm}(\omega)=\hat {\cal R}_{\pm j}(\omega) ={e^{{\omega
\over 2}}\over 2}{\sinh {j\omega \over 2} \over \cosh (k-1){\omega \over
2}\sinh{\omega \over 2}}, ~~j = 1, \ldots ,k-2 \non\\
&& \hat {\cal R}_{++}(\omega) = \hat {\cal R}_{--}(\omega)={e^{{\omega
\over 2}}\over 2}{\sinh {k\omega \over 2} \over 2\cosh{\omega \over 2}
\cosh (k-1){\omega \over 2}\sinh{\omega \over 2}}\non\\
&& \hat {\cal R}_{+-}(\omega) = \hat {\cal R}_{-+}(\omega)={e^{{\omega
\over 2}}\over 2}{\sinh (k-2){\omega \over 2} \over 2\cosh{\omega \over 2}
\cosh (k-1){\omega \over 2}\sinh{\omega \over 2}}
\label{RR12}
\end{eqnarray}

Let us now consider the particular solution D4 for the $so(4)$ case. The
corresponding Bethe Ansatz equations, as in the bulk, are basically two
copies of the $XXX$ spin chain equations ($so(4) = su(2) \otimes su(2)$),
namely
\begin{eqnarray}
e_{2\xi_{\tau}+1}(\lambda_{i}^{(\tau)})^{-1}
e_{1}(\lambda_{i}^{(\tau)})^{2N+1} &=& \prod_{ j=1}^{M^{(\tau)}}
e_{2}(\lambda_{i}^{(\tau)} - \lambda_{j}^{(\tau)})\
e_{2}(\lambda_{i}^{(\tau)} + \lambda_{j}^{(\tau)})
\label{d4ba}
\end{eqnarray}
where $\tau = \pm$. It is obvious that the only representations that remain
are the two two-dimensional spinor representations. The Bethe Ansatz
equations are then two decoupled equations, as it is also evident from
(\ref{RR1}): one has ${\cal R}_{+-} (\omega)= {\cal R}_{-+} (\omega)=0$;
moreover ${\cal R}_{++}(\omega)$, ${\cal R}_{--}(\omega)$ are given by
(\ref{RR12}) for $k=2$, and
\begin{eqnarray}
F^{\tau}(\lambda) = a_{1}(\lambda) +a_{2}(\lambda)+
\sum_{j=1}^{\nu^{(\tau)}}\Big(a_{2}(\lambda -\tilde \lambda_{j}^{(\tau)})+
a_{2}(\lambda +\tilde \lambda_{j}^{(\tau)})\Big), ~~ G^{\tau}(\lambda,
\xi_{\tau}) = -a_{2\xi_{\tau}+1}(\lambda). \label{xi4}
\end{eqnarray}

\subsection*{C. $ \bf sp(2k)$}
The ground state in this case consists of $k-1$ filled Dirac seas of
two-strings ($\lambda_{0}^{(j)} \pm{i \over 2}$), and the $k$ sea is filled
with real strings. The Bethe Ansatz equations take the form
\begin{eqnarray}
e_{2}(\lambda_{i}^{(1)})^{2N+1} &\!\!=\!\!& -\prod_{j=1}^{M^{(1)}}
e_{2}^{2}(\lambda_{i}^{(1)} - \lambda_{j}^{(1)})\
e_{2}^{2}(\lambda_{i}^{(1)} + \lambda_{j}^{(1)})\ e_{4}(\lambda_{i}^{(1)} -
\lambda_{j}^{(1)}) e_{4}(\lambda_{i}^{(1)} + \lambda_{j}^{(1)}) \non\\
&& \times \prod_{j=1}^{M^{(2)}}e_{-1}(\lambda_{i}^{(1)} -
\lambda_{j}^{(2)})\ e_{-1}(\lambda_{i}^{(1)} + \lambda_{j}^{(2)})\
e_{-3}(\lambda_{i}^{(1)} - \lambda_{j}^{(2)})\ e_{-3}(\lambda_{i}^{(1)} +
\lambda_{j}^{(2)})\,, \non\\
e_{2}(\lambda_{i}^{(l)}) &\!\!=\!\!& -\prod_{i j=1}^{M^{(l)}}
e_{2}^{2}(\lambda_{i}^{(l)} - \lambda_{j}^{(l)})\
e_{2}^{2}(\lambda_{i}^{(l)} + \lambda_{j}^{(l)})\ e_{4}(\lambda_{i}^{(l)} -
\lambda_{j}^{(l)}) e_{4}(\lambda_{i}^{(l)} + \lambda_{j}^{(l)}) \non\\
&& \times \prod_{\tau =\pm 1} \prod_{ j=1}^{M^{(l+\tau)}}
e_{-1}(\lambda_{i}^{(l)} - \lambda_{j}^{(l+\tau)})\
e_{-1}(\lambda_{i}^{(l)} + \lambda_{j}^{(l+\tau)}) e_{-3}(\lambda_{i}^{(l)}
- \lambda_{j}^{(l+\tau)})\ e_{-3}(\lambda_{i}^{(l)} +
\lambda_{j}^{(l+\tau)})\ \non\\
&& l= 2, \ldots , k-1 \non\\
e_{2}(\lambda_{i}^{(k)}) &\!\!=\!\!& -\prod_{ j=1}^{M^{(k)}}
e_{4}(\lambda_{i}^{(k)} - \lambda_{j}^{(k)})\
e_{4}(\lambda_{i}^{(k)} + \lambda_{j}^{(k)})\ \non\\
&& \times \prod_{j=1}^{M^{(k-1)}}e_{-1}(\lambda_{i}^{(k)} -
\lambda_{j}^{(k-1)})\ e_{-1}(\lambda_{i}^{(k)} + \lambda_{j}^{(k-1)})
e_{-3}(\lambda_{i}^{(k)} - \lambda_{j}^{(k-1)})\ e_{-3}(\lambda_{i}^{(k)} +
\lambda_{j}^{(k-1)}) \non\\
\label{BAE3d}
\end{eqnarray}
where for the general diagonal solutions of the reflection equation we have
to multiply:
\begin{itemize}
\item[--]
{\bf for D1}, the LHS of the $k^{th}$ equation with $-e_{2\xi+ \kappa}^{-1}$,
\item[--]
{\bf for D3}, the LHS of the $m^{th}$ equations with $e_{2\xi+ m
+1}^{-1}e_{2\xi+ m -1}^{-1}$, $m=1,\ldots ,k-1$.
\end{itemize}
Again the densities for the state with $\nu^{(l)}$ holes in the $l$ sea
satisfy the same equations (\ref{sigma}) as in the $so(2k+1)$ case with
$\sigma$ given by (\ref{cols}), $a^{j}(\lambda) = 2a_{2}(\lambda)
\delta_{j1}$ and
\begin{eqnarray}
\hat \epsilon^{j}(\omega) ={\cosh(k+1 -j){\omega \over 2}
\over 2\cosh{\omega \over 2} \cosh(k+1){\omega \over 2}},
~~j=1,\ldots,k \label{ener3}
\end{eqnarray}
\begin{eqnarray}
F^{j}(\lambda) &\!\!=\!\!& 3a_{2}(\lambda)+a_{4}(\lambda)
-2(a_{1}(\lambda)+a_{3}(\lambda))
+(a_{1}(\lambda)+a_{3}(\lambda))\delta_{j1} \non\\
&+& \sum_{j=1}^{\nu^{(l)}} \Big((2a_{2} +a_{4}) (\lambda -\tilde
\lambda_{j}^{(l)})+(2a_{2}+a_{4}) (\lambda +\tilde
\lambda_{j}^{(l)})\Big) \delta_{jl} \non\\
&-& \sum_{j=1}^{\nu^{(l)}} \Big((a_{1}+a_{3}) (\lambda -\tilde
\lambda_{j}^{(l)})+(a_{1}+a_{3}) (\lambda +\tilde \lambda_{j}^{(l)})\Big)
(\delta_{j,l+1}+\delta_{j,l-1}), \non\\
&& (j=1, \ldots ,k-1) \non \\
F^{k}(\lambda) &\!\!=\!\!& a_{2}(\lambda)+a_{4}(\lambda)
-(a_{1}(\lambda)+a_{3}(\lambda)) +\sum_{j=1}^{\nu^{(l)}} \Big(a_{4}(\lambda
-\tilde \lambda_{j}^{(l)})+a_{4}(\lambda +\tilde \lambda_{j}^{(l)})\Big)
\delta_{lk} \non\\
&-& \sum_{j=1}^{\nu^{(l)}} \Big((a_{1}+a_{3}) (\lambda -\tilde
\lambda_{j}^{(l)})+(a_{1}+a_{3}) (\lambda +\tilde \lambda_{j}^{(l)})\Big)
\delta_{l,k-1},
\label{f3}
\end{eqnarray}
\begin{eqnarray}
\mbox{\bf for D1:} && G^{j}(\lambda, \xi) = -a_{2\xi+\kappa
}(\lambda)\delta_{jk}, \non\\
\mbox{\bf for D3:} && G^{j}(\lambda, \xi) =
-(a_{2\xi+m+1}(\lambda)+a_{2\xi+m-1}(\lambda))\delta_{jm} ~~
(m=1,...,k-1)
\label{xi3}
\end{eqnarray}
The non-vanishing entries of $\hat {\cal K}$ are
\begin{eqnarray}
\hat {\cal K}_{ij}(\omega) &=& (1+\hat
a_{2}(\omega))^{2}\delta_{ij} - (\hat a_{1}(\omega)+\hat
a_{3}(\omega))(\delta_{i,j+1}+ \delta_{i,j-1}), ~~i,j =1,\ldots, k-1
\non\\ 
\hat {\cal K}_{k-1,k}(\omega) &=& \hat {\cal K}_{k,k-1}(\omega)
= - (\hat a_{1}(\omega)+\hat a_{3}(\omega)),
\non\\ 
\hat {\cal K}_{kk}(\omega) &=& 1+\hat a_{4}(\omega). 
\label{K3}
\end{eqnarray}
As in the previous cases the solution of (\ref{sigma}) has the form
(\ref{sigma2}), (\ref{ef1}) with
\begin{equation}
\hat {\cal R}_{ij}(\omega) ={ e^{\omega} \over 2\cosh {\omega \over 2}}
{\sinh \min(i,j) {\omega \over 2} \cosh(k+1 -\max(i,j)){\omega \over 2}
\over \cosh (k+1){\omega \over 2}\sinh{\omega \over 2} }, ~~i,j = 1, \ldots
,k
\label{RR3}
\end{equation}

\section{Scattering}
\setcounter{equation}{0}

Having obtained the excitations of the model, we are ready to compute the
complete boundary $S$ matrix. For this purpose we follow the formulation
developed by Korepin, and later by Andrei and Destri \cite{K,AD}. First we
have to implement the so-called quantisation condition,
\begin{equation}
(e^{2iNp^{l}}S-1)|\tilde \lambda_{i}^{(l)} \rangle = 0
\label{qc1}
\end{equation}
where $p^{l}$ is the momentum of the particle (in our case, the hole) with
rapidity $\tilde \lambda_{1}^{l}$. For the case of $\nu$ (even) holes in
$l$ sea we insert the integrated density (\ref{sigma2}) into the
quantisation condition (\ref{qc1}). We use the dispersion relation
\begin{equation}
\epsilon^{l}(\lambda) = {1 \over  2\pi} {d \over d \lambda} \
p^{l}(\lambda)
\end{equation}
and the sum rule $\displaystyle N\int_{0}^{\tilde \lambda_{i}} d\lambda
\sigma(\lambda) \in \ZZ_{+}$. The boundary $S$ matrix for the right
boundary is taken to be one of the diagonal solutions D1, D2, D3, D4,
whereas the boundary $S$ matrix for the left boundary is proportional to
unit. We end up with the following expression for the boundary scattering
amplitudes:
\begin{eqnarray}
\alpha^{+l}
\alpha^{-l} = \exp \Bigl \{ 2 \pi N \int_{0}^{\tilde
\lambda_{1}}d\lambda \Bigl (\sigma^{l}(\lambda)
-2\epsilon^{l}(\lambda)\Bigr) \Bigr \}
\end{eqnarray}
with
\begin{eqnarray}
\alpha^{-l}(\lambda, \xi) = k_{0}(\lambda) k_{1}(\lambda, \xi),
 ~~\alpha^{+l}(\lambda) = k_{0}(\lambda)
\label{sol1}
\end{eqnarray}
where $\alpha^{-l}$ is the first element of the diagonal boundary $S$
matrix. It is obvious that $\alpha^{+l}$ has no $\xi$ dependence and
realises just the overall factor in front of the unit matrix at the left
boundary (recall that $K^{+} =1$). Moreover,
\begin{eqnarray}
k_{0}(\tilde \lambda_{1}^{l}) = \exp \Big \{ i\pi
\int_{0}^{\tilde \lambda_{1}^{l}} d\lambda \Phi_{0}^{l}(\lambda) \Big \},
~~k_{1}(\tilde \lambda_{1}^{l}, \xi) = \exp \Big \{ 2i\pi \int_{0}^{\tilde
\lambda_{1}^{l}} d\lambda \Phi_{1}^{l}(\tilde \lambda_{1}^{l}, \xi)\Big \}
\label{ampl}
\end{eqnarray}
with $\Phi_{0,1}^{l}$ given by (\ref{ef1}). We finally restrict ourselves
to $l=1$ in the first sea and we write the latter expression in term of the
Fourier transform (\ref{ef1}) of $\Phi_{0,1}^{1}$,
\begin{eqnarray}
k_{0}(\lambda) = \exp\Big \{ -{1\over 2}
\int_{-\infty}^{\infty} {d\omega \over \omega} \hat
\Phi_{0}^{1}(\omega)e^{-i\omega \lambda} \Big \}, ~~k_{1}(\lambda,
\xi) = \exp \Big \{ -\int_{-\infty}^{\infty} {d\omega \over
\omega} \hat \Phi_{1}^{1}(\omega, \xi)e^{-i\omega \lambda} \Big
\}. \label{ampl2}
\end{eqnarray}
In what follows we express the scattering amplitudes in terms of
$\Gamma$-functions.

\subsection*{A. $\bf so(n)$}
Before we write down the explicit expressions for the boundary $S$
matrices, let us recall the form of the exact bulk $S$ matrix. It is easy
to compute the scattering amplitude between two holes (vectors) in the
first sea. The bulk scattering amplitude comes from the contribution of the
terms of $\Phi_{0}^{1}$ given by eqs. (\ref{f1}), (\ref{ef1}), (\ref{RR1}),
(\ref{f2}), with argument $\lambda \pm \tilde \lambda_{j}$. After some
algebra and using the following identity
\begin{eqnarray}
{1 \over 2}\int_{0}^{\infty}{d\omega \over \omega}{e^{-{\mu\omega\over 2}}
\over \cosh{\omega \over 2}} = \ln \, {\Gamma({\mu+1\over 4}) \over
\Gamma({\mu+3\over 4})}
\label{gamma1}
\end{eqnarray}
we conclude that the hole--hole scattering amplitude is given by the
expression
\begin{eqnarray}
S_{0}(\lambda) = {\tan \pi ({i\lambda -1\over n-2}) \over \tan \pi
({i\lambda +1\over n-2})} \; {\Gamma({i\lambda \over n-2}) \over
\Gamma({-i\lambda \over n-2})} \; {\Gamma({-i\lambda \over n-2}+{1\over 2})
\over \Gamma({i\lambda \over n-2}+{1\over 2})} \; {\Gamma({-i\lambda +
1\over n-2}) \over \Gamma({i\lambda +1\over n-2})} \; {\Gamma({i\lambda +1
\over n-2}+{1\over 2}) \over \Gamma({-i\lambda +1 \over n-2}+{1\over 2})}.
\label{bulk1}
\end{eqnarray}

The explicit bulk $S$ matrix then has the following structure \footnote{We
do not compute all the eigenvalues of the bulk $S$ matrix via the Bethe
Ansatz. Such a general computation, in the bulk, is beyond the scope of
this work.}
\begin{eqnarray}
S(\lambda) = {S_{0}(\lambda) \over
(i\lambda+\kappa)(i\lambda+1)} (i\lambda(i\lambda+\kappa)\II +(i\lambda
+\kappa) P - i\lambda Q )
\label{smatrix}
\end{eqnarray}
Now, we give the expressions for the boundary $S$ matrix, which follow from
(\ref{ampl}), (\ref{ampl2}), and the duplication formula for the $\Gamma$
function
\begin{eqnarray}
2^{2x-1}\Gamma(x+{1 \over 2})\Gamma(x)= \pi^{{1\over 2}}\Gamma(2x).
\label{gamma3}
\end{eqnarray}
The $\xi$-independent part of the overall factor $k_{0}$, eq. (\ref{ampl}),
is given by
\begin{eqnarray}
k_{0}(\lambda) = Y_{0}(\lambda)  {\Gamma({i\lambda
\over n-2}) \over \Gamma({-i\lambda \over n-2})} \; {\Gamma({-i\lambda
\over n-2}+{3\over 4}) \over \Gamma({i\lambda \over n-2}+{3\over
4})} \; {\Gamma({i\lambda +1/2\over n-2}+{3\over 4})
\over \Gamma({-i\lambda+1/2\over n-2} + {3\over
4})} \; {\Gamma({-i\lambda+1/2\over n-2} +{1\over 2})
\over \Gamma({i\lambda +1/2\over n-2}+{1\over 2})}
\label{bound1}
\end{eqnarray}
where
\begin{eqnarray}
Y_{0}(\lambda) = {\sin \pi ({i\lambda+1/2 \over n-2} -{1\over
4}) \over \sin \pi ({i\lambda-1/2\over n-2} +{1\over 4})} \;
{\sin \pi ({i\lambda -1/2\over n-2} +{1\over 2}) \over \sin
\pi ({i\lambda+1/2 \over n-2} -{1\over 2})} \; {\sin \pi
({i\lambda \over n-2}+{1\over 4}) \over \sin \pi ({i\lambda \over
n-2}-{1\over 4})} \,. \label{Y00}
\end{eqnarray}
Note that our solution includes the necessary CDD factors both for the bulk
and boundary matrices (see also \cite{McS,ogiev}). The expression for the
$\xi$-dependant part $k_{1}$ depends on which solutions D1, D2, D3, D4 we
consider:

\bigskip

{\bf D1:}
\begin{eqnarray}
k_{1}(\lambda,\xi')= {\Gamma({i\lambda +\xi'\over n-2} +{1\over
2}) \over \Gamma({-i\lambda+\xi' \over n-2}+{1\over 2})} \;
{\Gamma({-i\lambda+\xi' \over n-2}+1) \over \Gamma({\xi'+i\lambda
\over n-2}+1 )}
\label{D1}
\end{eqnarray}
where $\xi' =\xi -{1\over 2}$ is the renormalised boundary parameter (see
also \cite{done}). We also compute the $\beta$ element of the $K$ matrix
(\ref{1}) by employing the ``duality'' transformation $\xi \to -\xi$, and
the symmetry of the $K$ matrix and of the transfer matrix (see also
\cite{done}). In particular, we obtain a set of Bethe Ansatz equations for
$\xi \to -\xi$, which allows to determine the difference
\begin{eqnarray}
\hat \Phi_{1}^{1}(\omega, -\xi) - \hat \Phi_{1}^{1}(\omega, \xi) =
e^{-(2\xi -1){\omega \over 2}}
\end{eqnarray}
and consequently
\begin{eqnarray}
{\beta^{-}(\lambda, \xi) \over \alpha^{-}(\lambda, \xi)} =
e_{2\xi'}(\lambda).
\end{eqnarray}
This provides us with a consistency check for our procedure since one
obtains independently the exact ratio between different elements of the
reflection matrix. Thus, we have completely determined the $K$ matrix that
corresponds to the solution D1.

\bigskip

{\bf D2:}
\begin{eqnarray}
k_{1}(\lambda,\xi')&=& {\tan \pi ({i\lambda-\xi_{n}' \over n-2})
 \over \tan \pi ({i\lambda+\xi_{n}' \over n-2})}  \;
{\Gamma({i\lambda+\xi_{1}' \over n-2}+{1\over 2} ) \over
\Gamma({-i\lambda+\xi_{1}' \over n-2}+{1\over 2})} \;
{\Gamma({-i\lambda+\xi_{1}' \over n-2}+1) \over
\Gamma({i\lambda+\xi_{1}' \over n-2}+1 )}
 {\Gamma({i\lambda +\xi_{n}'\over n-2} +{1\over 2})
\over \Gamma({-i\lambda +\xi_{n}'\over n-2}+{1\over 2})} \;
{\Gamma({-i\lambda +\xi_{n}'\over n-2}) \over \Gamma({i\lambda
+\xi_{n}'\over n-2})}
\label{D2}
\end{eqnarray}
where $\xi_{1}' =\xi_{1} -{1\over 2}$ and $\xi_{n}' =\xi_{n} +{1\over 2}$
are the renormalised boundary parameters. Remark that the constraint
$\xi_{1} +\xi_{n} = \kappa-1$ is also true for the renormalised boundary
parameters, namely $\xi_{1}' + \xi_{n}' = \kappa -1$. Similarly we employ
the duality transformation which, for this solution, reads $\xi_{1} \to
\xi_{n}$. We then determine the difference
\begin{eqnarray}
\hat \Phi_{1}^{1}(\omega, \xi_{n}, \xi_{1}) - \hat
\Phi_{1}^{1}(\omega, \xi_{1}, \xi_{n}) = e^{-(2\xi_{1} -1){\omega
\over 2}}+e^{(2\xi_{n} +1){\omega \over 2}}
\end{eqnarray}
and the ratio
\begin{eqnarray}
{\gamma^{-}(\lambda, \xi_{1}, \xi_{n}) \over
\alpha^{-}(\lambda, \xi_{1}, \xi_{n})} = e_{2
\xi_{1}'}(\lambda)e_{-2 \xi_{n}'}(\lambda).
\end{eqnarray}
Again we have here an independent consistency check. \\
Of course, we need to determine one further ratio, namely
${\beta^{-}(\lambda, \xi_{1}, \xi_{n}) \over \alpha^{-}(\lambda, \xi_{1},
\xi_{n})}$. We have not been able to explicitly extract this information
from the Bethe Ansatz formulation. Nevertheless, since the $K$ matrix is a
solution of the reflection equation, the only choice we have for this ratio
is
\begin{eqnarray}
{\beta^{-}(\lambda, \xi_{1},
\xi_{n}) \over \alpha^{-}(\lambda, \xi_{1}, \xi_{n})} = e_{2
\xi_{1}'}(\lambda).
\end{eqnarray}
With this, we have completed the derivation of the $K$ matrix that
corresponds to the solution D2.

\bigskip

{\bf D3:} The $\xi$-dependent part of the $K$ matrix overall factor that
corresponds to D3 is
\begin{eqnarray}
k_{1}(\lambda,\xi') = {\Gamma({i\lambda +\xi'\over n-2}
+{1\over 2}) \over \Gamma({-i\lambda +\xi' \over n-2}
+{1\over 2})} \; {\Gamma({-i\lambda +\xi'\over n-2}+1)
\over \Gamma({i\lambda +\xi'\over n-2}+1 )} \;
{\Gamma({-i\lambda+1/2 \over n-2}+{3\over 4}) \over
\Gamma({i\lambda+1/2\over n-2}  + {3\over 4})} \;
{\Gamma({i\lambda+1/2\over n-2} +{1\over 4}) \over
\Gamma({-i\lambda+1/2 \over n-2}+{1\over 4})} \,.
\label{D3}
\end{eqnarray}
The total overall factor for this solution is
\begin{eqnarray}
k_{0}(\lambda)k_{1}(\lambda,\xi') &=& Y_{1}(\lambda) {\Gamma({i\lambda
\over n-2}) \over \Gamma({-i\lambda \over n-2})} \; {\Gamma({-i\lambda
\over n-2}+{3\over 4}) \over \Gamma({i\lambda \over n-2}+{3\over 4})} \;
{\Gamma({i\lambda +1/2\over n-2}+ {1\over 4}) \over
\Gamma({-i\lambda+1/2\over n-2} + {1\over 4})} \;
{\Gamma({-i\lambda +1/2\over n-2} +{1\over 2}) \over
\Gamma({i\lambda  +1/2\over n-2}+{1\over 2})} \non\\
&& \times {\Gamma({i\lambda +\xi'\over n-2} +{1\over 2}) \over
\Gamma({-i\lambda +\xi' \over n-2}+{1\over 2})} \;
{\Gamma({-i\lambda+\xi' \over n-2}+1) \over \Gamma({i\lambda
+\xi'\over n-2}+1)}
\label{all1}
\end{eqnarray}
where
\begin{eqnarray}
Y_{1}(\lambda) = {\sin \pi ({i\lambda +1/2\over n-2} -{1\over
4}) \over \sin \pi ({i\lambda-1/2\over n-2} +{1\over 4})} \;
{\sin \pi ({i\lambda -1/2\over n-2} +{1\over 2}) \over \sin
\pi ({i\lambda +1/2\over n-2} -{1\over 2})} \; {\sin \pi
({i\lambda \over n-2}+{1\over 4}) \over \sin \pi ({i\lambda \over
n-2}-{1\over 4})} \,.
\label{Y10}
\end{eqnarray}
The latter expression agrees with the one found in \cite{McS}, although we
find slightly different CDD factors. Again $\xi' =\xi -{1\over 2}$ is the
renormalised boundary parameter which now must have the fixed value
$\xi'={n\over 4} -m$, to ensure the resulting $K$ matrix satisfies the
reflection equation. Unfortunately, in this case, we cannot employ a
duality transformation to derive the ratio ${\beta^{-} \over
\alpha^{-}}$ (\ref{3}). However, since the $K$ matrix is a solution of the
reflection equation the ratio must be
\begin{eqnarray}
{\beta^{-}(\lambda) \over \alpha^{-}(\lambda)}= e_{2
\xi'}(\lambda).
\end{eqnarray}

\bigskip

{\bf D4:} For this solution the $\xi$-independent part is given by
(\ref{bound1}) with $n=4$, while the $\xi$-dependent part is given by
($\tau=\pm$):
\begin{eqnarray}
k_{1\tau}(\lambda,\xi'_{\tau})= {\Gamma({i\lambda+\xi_{\tau}' \over 2}
+{1\over 4}) \over \Gamma({-i\lambda+\xi_{\tau}' \over 2}+{1\over 4})} \;
{\Gamma({-i\lambda+\xi_{\tau}' \over 2}+{3\over 4}) \over
\Gamma({i\lambda +\xi_{\tau}'\over 2}+{3\over 4} )}. \label{D4}
\end{eqnarray}
We exploit the duality transformation one more time, $\xi_{\tau} \to -
\xi_{\tau}$, to calculate the ratio
\begin{eqnarray}
{\beta_{\tau}^{-}(\lambda) \over \alpha_{\tau}^{-}(\lambda)}= e_{2
\xi_{\tau}'}(\lambda)
\end{eqnarray}
with the renormalised boundary parameter $\xi_{\tau}'=\xi_{\tau}-{1\over
2}$. We recall that the two spinor representations are two-dimensional, and
the corresponding $K$ matrices are of course two-dimensional with
\begin{eqnarray}
K_{\tau}^{-}(\lambda, \xi_{\tau'}) =
diag(\alpha_{\tau}^{-}(\lambda,
\xi_{\tau}'),\beta_{\tau}^{-}(\lambda, \xi_{\tau}')) \,.
\end{eqnarray}
In other words, we have obtained two copies of the $XXX$ boundary $S$
matrix, with two free boundary parameters $\xi_{\tau}'$.

\subsection*{B. $\bf sp(n)$}
The corresponding bulk scattering amplitude for the $sp(n)$ is given by the
following expression
\begin{eqnarray}
S_{0}(\lambda) = {\tan \pi ({i\lambda -1\over n+2}) \over \tan \pi
({i\lambda +1\over n+2})} \; {\Gamma({i\lambda \over n+2}) \over
\Gamma({-i\lambda \over n+2})} \; {\Gamma({-i\lambda \over n+2}+{1\over 2})
\over \Gamma({i\lambda \over n+2}+{1\over 2})} \; {\Gamma({-i\lambda +
1\over n+2}) \over \Gamma({i\lambda +1\over n+2})} \; {\Gamma({i\lambda +1
\over n+2}+{1\over 2}) \over \Gamma({-i\lambda +1 \over n+2}+{1\over 2})}
\; S_b(\lambda)
\label{bulk2}
\end{eqnarray}
where
\begin{equation}
  \label{eq:Sblambda}
  S_b(\lambda) = \exp\left[ -\int_{-\infty}^{\infty}
  \frac{d\omega}{\omega} 
  \frac{\cosh \frac{k \omega}{2} \; e^{-i\omega \lambda}}
  {2 \cosh \frac{\omega}{2} \cosh\frac{(k+1)\omega}{2}}
  \right]
\end{equation}
Again the structure of the $S$ matrix is given by (\ref{smatrix}) with
$S_{0}$ given above. 
$k_0(\lambda)$ is again given by the integral representation
(\ref{ampl2}) using (\ref{ef1}), (\ref{f3}) and (\ref{RR3}).
The corresponding $\xi$-dependent parts for each solution are given
respectively by:

\bigskip

{\bf D1:}
\begin{eqnarray}
k_{1}(\lambda,\xi')= {\Gamma({i\lambda+\xi' \over n+2} +{1\over
2}) \over \Gamma({-i\lambda +\xi'\over n+2}+{1\over 2})} \;
{\Gamma({-i\lambda +\xi'\over n+2}+1) \over \Gamma({i\lambda
+\xi'\over n+2}+1 )}
\label{D12}
\end{eqnarray}
where $\xi' =\xi -1$ is the renormalised boundary parameter. We also
compute the $\beta$ element of the $K$ matrix by employing the duality
transformation ($\xi \to -\xi$) and we obtain the difference
\begin{eqnarray}
\hat \Phi_{1}^{1}(\omega, -\xi) - \hat
\Phi_{1}^{1}(\omega, \xi) = e^{-(2\xi -2){\omega \over 2}}
\end{eqnarray}
and consequently
\begin{eqnarray}
{\beta^{-}(\lambda, \xi) \over
\alpha^{-}(\lambda, \xi)} = e_{2 \xi'}(\lambda).
\end{eqnarray}

\bigskip

{\bf D3:}
\begin{eqnarray}
k_{1}(\lambda,\xi') &=& {\tan \pi({i\lambda -1/2\over n+2 }
-{1\over 4}) \over \tan \pi({i\lambda +1/2\over n+2 } +{1\over
4})} \; {\Gamma({i\lambda +\xi'\over n+2} +{1\over 2}) \over
\Gamma({-i\lambda +\xi'\over n+2}+{1\over 2})} \;
{\Gamma({-i\lambda +\xi'\over n+2}+1) \over \Gamma({i\lambda
+\xi'\over n+2}+1 )} 
{\Gamma({-i\lambda+1/2 \over n+2}+{1\over 4}) \over
\Gamma({i\lambda+1/2\over n+2} + {1\over 4})} \;
{\Gamma({i\lambda+1/2\over n+2} +{3\over 4}) \over
\Gamma({-i\lambda+1/2 \over n+2}+{3\over 4})} \,.\qquad
\label{D32}
\end{eqnarray}
$\xi' =\xi-1$ is the renormalised boundary parameter which must have the
fixed value $\xi'={n\over 4} -m$, in order for the resulting $K$ matrix to
satisfy the reflection equation. This result also agrees with 
the $\xi$-dependent part of the corresponding result of \cite{McS}.
Unfortunately we cannot employ a
duality symmetry in this case to derive the ratio ${\beta \over
\alpha}$ (\ref{3}), but again, since the $K$ matrix is a solution of the
reflection equation, the ratio must be
\begin{eqnarray}
{\beta^{-}(\lambda) \over \alpha^{-}(\lambda)}= e_{2 \xi'}(\lambda) \,.
\end{eqnarray}

\section{Conclusion}
\setcounter{equation}{0}

We have established here the classification of diagonal, antidiagonal
and mixed solutions to the reflection
equations based on the Yangian $R$ matrices for Lie (super)algebras
$so(m)$, $sp(n)$ and $osp (m|n)$.  The next step would then be to
obtain the complete classification of $K$ matrices, still associated
with one-dimensional quantum boundary space, which is technically
more complicated.
The cases of higher dimensional bulk representations, as well as
operator valued $K$-matrices remain to be classified.

Some of the solutions we obtained have then been used to compute the
spectrum and 
scattering datas for integrable spin chain systems with non trivial
boundary conditions. Only the Lie algebra case has been considered here;
the case of super Lie algebras $osp (m|n)$ can now be envisioned:
indeed, a suitable redefinition of indices allows us to consider a
similar exact ferromagnetic pseudo-vacuum (see \cite{martins}),
providing us with the starting point for the analytical Bethe Ansatz. 

Remark also that the analytical Bethe Ansatz construction proposed
here remains valid for upper triangular reflection matrices. Such
matrices can for instance be obtained by conjugating our diagonal
solutions by triangular matrices (following lemma \ref{lemma:invar}).

In addition, the use of the analytical Bethe Ansatz method precluded
applications to non-diagonal reflection matrices, for which no exact
pseudo-vacuum state is available. However, this situation
can be approached using the methods developed e.g. in \cite{bazh,nepo2002}.
For an alternative approach, see \cite{CLSW}.

\bigskip

\textbf{Acknowledgements:} We warmfully thank E.~Sokatchev for
discussions on symmetric spaces. 
This work was supported by the TMR Network
``EUCLID. Integrable models and applications: from strings to condensed
matter'', contract number HPRN-CT-2002-00325. J.A. wishes to thank LAPTH
for kind hospitality.


\begin{thebibliography}{99}

\bibitem{Che}
  I.V.~Cherednik, \textsl{Factorizing particles on a half line and root
  systems,} Theor. Math. Phys. \textbf{61} (1984) 977.

\bibitem{Skl}
  E.K.~Sklyanin, \textsl{Boundary conditions for integrable quantum
  systems,} J. Phys. \textbf{A21} (1988) 2375.

\bibitem{GhZ}
  S.~Ghoshal and A.B.~Zamolodchikov, \textsl{Boundary $S$ matrix and
  boundary state in two-dimensional integrable quantum field theory,} Int.
  Journ. Mod. Phys. \textbf{A9} (1994) 3841 and \texttt{hep-th/9306002}.

\bibitem{dVG}
  H.J.~de~Vega and A.~Gonz\'{a}lez-Ruiz, \textsl{Boundary $K$ matrices for
  the six vertex and the $N(2N-1)$ $A(N-1)$ vertex models,} J. Phys.
  \textbf{A26} (1993) L519 and \texttt{hep-th/9211114}.

\bibitem{ZF}
  A.~B.~Zamolodchikov and A.~B.~Zamolodchikov, \textsl{Factorized
  S-matrices in two dimensions as the exact solutions of certain
  relativistic quantum field theory models,} Annals Phys.{\bf 120} (1979)
  253; L.~D.~Faddeev, \textsl{Quantum completely integrable models in field
  theory,} Sov. Sci. Rev. \textbf{C1} (1980) 107.

\bibitem{Mintchev}
  A.~Liguori, M.~Mintchev and L.~Zhao, \textsl{Boundary exchange algebras
  and scattering on the half-line,} Commun. Math. Phys.\textbf{194} (1998)
  569 and \texttt{hep-th/9607085}.

\bibitem{DoMu}
  J.~Donin and A.I.~Mudrov, \textsl{Reflection equation, twist, and
  equivariant quantization,} \texttt{math.QA/0204295}.

\bibitem{DoKuMu}
  J.~Donin, P.P.~Kulish and A.I.~Mudrov, \textsl{On universal solution to
  reflection equation,} \texttt{math.QA/0210242}.

\bibitem{BaKo}
  P.~Baseilhac and K. Koizumi, \textsl{Sine--Gordon quantum field
  theory on the half line with quantum boundary degrees of freedom,}
  Nucl. Phys. \textbf{B649} (2003) 491 and \texttt{hep-th/0208005}.

\bibitem{Ma}
  G.W.~Delius and N.J.~MacKay, \textsl{Quantum group symmetry in
  sine-Gordon and affine Toda field theories on the half-line,} Commun.
  Math. Phys. \textbf{233} (2003) 173 and \texttt{hep-th/0112023}.

\bibitem{DeNe}
  G.W.~Delius and R.I.~Nepomechie, \textsl{Solutions of the boundary
    Yang--Baxter equation for arbitrary spin,} J. Phys. \textbf{A35} (2002)
  L341 and \texttt{hep-th/0204076}.

\bibitem{Ne}
  R.I.~Nepomechie, \textsl{Boundary quantum group generators of type A,}
  Lett. Math. Phys. \textbf{62} (2002) 83 and \texttt{hep-th/0204181}.

\bibitem{AR}
  J.~Abad and M.~Rios, \textsl{Nondiagonal solutions to reflection
    equations in $SU(N)$ spin chains,} Phys. Lett. \textbf{B352} (1995)
  92 and \texttt{hep-th/9502129}.

\bibitem{MeNe}
  L.~Mezincescu and R.I.~Nepomechie, \textsl{Fractional-spin integrals
  of motion for the boundary Sine--Gordon model at the free fermion
  point,} Int. J. Mod. Phys. \textbf{A13} (1998) 2747 and
  \texttt{hep-th/9709078}. 

\bibitem{Li}
  A.~Lima-Santos, \textsl{Reflection $K$ Matrices for 19-vertex models,}
  Nucl. Phys. \textbf{B558} (1999) 637 and \texttt{solv-int/9906003}.

\bibitem{Ga}
  G.M.~Gandenberger, \textsl{New non-diagonal solutions to the $a_{n}^{(1)}$
    boundary Yang--Baxter equation,} \texttt{hep-th/9911178}.

\bibitem{Kim}
  J.D. Kim, \textsl{Boundary $K$ matrix for the quantum
  Mikhailov-Shabat model,} \texttt{hep-th/9412192}.

\bibitem{Dr}
  V.G.~Drinfel'd, \textsl{Hopf algebras and the quantum Yang--Baxter
  equation,} Soviet. Math. Dokl. \textbf{32} (1985) 254 and \textsl{A new
  realization of {Y}angians and quantized affine algebras,} Soviet. Math.
  Dokl. \textbf{36} (1988) 212.

\bibitem{soya}
  D.~Arnaudon, J.~Avan, N.~Cramp\'{e}, L.~Frappat and {\'E}.~Ragoucy,
  \textsl{$R$ matrix presentation for (super)-Yangians $Y(g)$,} J. Math.
  Phys. \textbf{44} (2003) 302 and \texttt{math.QA/0111325}.

\bibitem{LAPTH844}
  M.~Mintchev, {\'E}.~Ragoucy and P.~Sorba, \textsl{Spontaneous symmetry
  breaking in the gl({N})-{NLS} hierarchy on the half line,} J. Phys.
  \textbf{A34} (2001) 8345 and {\tt hep-th/0104079}.

\bibitem{McS}
  N.J.~MacKay and B.J.~Short, \textsl{Boundary scattering, symmetric spaces
  and the principal chiral model on the half-line,} Commun. Math. Phys.
  \textbf{233} (2003) 313 and \texttt{hep-th/0104212}.

\bibitem{Mori}
  M.~Moriconi, \textsl{Integrable boundary conditions and reflection
  matrices for the $O(N)$ nonlinear sigma model,} Nucl. Phys. \textbf{B619}
  (2001) 396 and \texttt{hep-th/0108039}.

\bibitem{olsh}
  G.I.~Olshanski, \textsl{Twisted Yangians and infinite dimensional
    Lie algebras}, in ``Quantum groups'', Lecture Notes in Math. \textbf{1510}
  (P. Kulish ed.), pp. 104, NY 1992.

\bibitem{MNO}
  A.~Molev, M.~Nazarov and G.~Olshanski, \textsl{Yangians and classical
  {L}ie algebras,} Russ. Math. Surveys \textbf{51} (1996) 205 and
  \texttt{hep-th/9409025}.

\bibitem{Nou}
  M.~Noumi, \textsl{Macdonald's symmetric polynomials as zonal spherical
  functions on some quantum homogeneous spaces,} Adv. Math. \textbf{123}
  (1996) 16.

\bibitem{MRS}
  A.~Molev, {\'E}.~Ragoucy and P.~Sorba, \textsl{Coideal subalgebras in
  quantum affine algebras,} \texttt{math.QA/0208140}, Rev. Math. Phys,
  to appear.

\bibitem{ogiev}
  E.~Ogievetsky, N.Yu.~Reshetikhin and P.~Wiegmann, \textsl{The principal
  chiral field in two-dimensions on classical Lie algebras: the Bethe
  Ansatz solution and factorized theory of scattering,} Nucl. Phys.
\textbf{B280} (1987) 45.

\bibitem{BFKZ}
  M.T.~Batchelor, V.~Fridkin, A.~Kuniba and Y.K.~Zhou,
  \textsl{Solutions of the reflection equation for face and vertex
  models associated with $A_{n}^{(1)}$, $B_{n}^{(1)}$, $C_{n}^{(1)}$,
  $D_{n}^{(1)}$ and $A_{n}^{(2)}$,} Phys. Lett. \textbf{B376} (1996)
266 and \texttt{hep-th/9601051}.

\bibitem{Zirn} M. Zirnbauer, {\sl Riemannian symmetric superspaces and their
origin in random-matrix theory}, J. Math. Phys. {\bf 37} (1996) 4986
and
{\tt math-ph/9808012}.

\bibitem{KS}
  P.P.~Kulish and E.K.~Sklyanin, \textsl{Solutions of the Yang--Baxter
  equation,} Zap. Nauchn. Sem. LOMI, \textbf{95} (1980) 129 and J. Sov.
  Math. \textbf{19} (1982) 1596.

\bibitem{reshe}
  V.I.~Vichirko and N.Yu.~Reshetikhin, \textsl{Excitation spectrum of
  the anisotropic generalization of an $SU_3$ magnet,}
  Theor. Math. Phys. \textbf{56}  (1983) 805;\\
 N.Yu. Reshetikhin, \textsl{A method of functional equations
  in the theory of exactly solvable quantum systems,} Lett. Math. Phys.
  \textbf{7} (1983) 205; 
  \textsl{Integrable models of quantum one-dimensional magnets with $O(n)$
    and $Sp(2k)$ symmetry,} Theor. Math. Phys. \textbf{63} (1985) 555;
  \textsl{The spectrum of the transfer matrices connected with
    Kac--Moody algebras,} Lett. Math. Phys. \textbf{14} (1987) 235.

\bibitem{mnanal}
  L.~Mezincescu and R.I.~Nepomechie, \textsl{Analytical Bethe Ansatz for
  quantum algebra invariant spin chains,} Nucl. Phys. \textbf{B372} (1992)
  597.

\bibitem{KunSuz} 
  A. Kuniba and J. Suzuki, \textsl{Analytic Bethe Ansatz for
    fundamental representations of Yangians,}
  Commun. Math. Phys. \textbf{173} (1995) 225
  and \texttt{hep-th/9406180}.

\bibitem{AMN1}
  S. Artz, L. Mezincescu and R.I.~Nepomechie, \textsl{Spectrum of
  transfer matrix for $U_q(B_n)$ invariant $A_{2n}^{(2)}$ open spin
  chains,} Int. J. Mod. Phys. \textbf{A10} (1995) 1937, and
\texttt{hep-th/9409130}. 

\bibitem{AMN2}
  S. Artz, L. Mezincescu and R.I.~Nepomechie, \textsl{Analytical Bethe
    Ansatz for $A_{2n-1}^{(2)}$, $B_{n}^{(1)}$, $C_{n}^{(1)}$,
    $D_{n}^{(1)}$ quantum algebra invariant open spin chains,}
  J.~Phys. \textbf{A28} (1995) 5131 and
  \texttt{hep-th/9504085}. 

\bibitem{doikou}
  A.~Doikou, \textsl{Fusion and analytical Bethe Ansatz for the
  $A_{n-1}^{(1)}$ open spin chain,} J. Phys. \textbf{A33} (2000) 4755;
  \textsl{Quantum spin chain with ``soliton nonpreserving'' boundary
  conditions,} J. Phys. \textbf{A33} (2000) 8797.

\bibitem{baxter}
  R.J.~Baxter, \textsl{Partition function of the eight-vertex lattice
  model,} Ann. Phys. \textbf{70} (1972) 193; J. Stat. Phys. \textbf{8}
  (1973) 25; {\it Exactly solved models in statistical mechanics} (Academic
  Press, 1982)

\bibitem{korepin}
  V.E.~Korepin, \textsl{New effects in the massive Thirring model: 
repulsive case}, Comm. Math. Phys. \textbf{76} (1980) 165; 
V.E. Korepin, G.
  Izergin and N.M. Bogoliubov, {\it Quantum inverse scattering method,
  correlation functions and algebraic Bethe Ansatz} (Cambridge University
  Press, 1993).

\bibitem{mnfusion}
  L.~Mezincescu and R.I.~Nepomechie, \textsl{ Fusion procedure for open
  chains,} J. Phys. \textbf{A25} (1992) 2533.

\bibitem{Zhou}
  Y-K. Zhou, \textsl{Row transfer matrix functional relations for
  Baxter's eight-vertex and six-vertex models with open boundaries via
  more general reflection matrices,} 
Nucl. Phys. \textbf{B458} (1996) 504, \texttt{hep-th/9510095}.

\bibitem{done}
  A.~Doikou and R.I.~Nepomechie, \textsl{Bulk and boundary S matrices for
  the su(N) chain,} Nucl. Phys. \textbf{B521} (1998) 547 and
  \texttt{hep-th/9803118}; \textsl{Duality and quantum algebra symmetry of
  the $A_{n-1}^{(1)}$ open spin chain with diagonal boundary fields,}
  Nucl. Phys. \textbf{B530} (1998) 641 and \texttt{hep-th/9807065}.

\bibitem{GMN}
   M.T. Grisaru, L. Mezincescu, R.I. Nepomechie, \textsl{Direct
   Calculation of the Boundary $S$ Matrix for the Open Heisenberg
   Chain,} J. Phys. \textbf{A28} (1995) 1027. 

\bibitem{K}
  V.~Korepin, \textsl{Direct calculation of the $S$ matrix in the massive
  Thirring model,} Theor. Math. Phys. \textbf{41} (1979) 953.

\bibitem{AD}
  N.~Andrei and C.~Destri, \textsl{Dynamical symmetry breaking and
  fractionization in a new integrable model,} Nucl. Phys. \textbf{B231}
  (1984) 445.

\bibitem{martins} 
  M.J.~Martins, \textsl{Bethe Ansatz solution of the $Osp(1|2n)$
  invariant spin chain,} Phys. Lett. \textbf{B359} (1995) 334.

\bibitem{bazh}
  V.V.~Bazhanov and N.Yu.~Reshetikhin, \textsl{Critical RSOS models and
  conformal field theory}, Int. J. Mod. Phys. \textbf{A4} (1989) 115.

\bibitem{nepo2002}
  R.I.~Nepomechie, \textsl{Functional relations and Bethe Ansatz for the
  $XXZ$ chain,} \texttt{hep-th/0211001}.

\bibitem{CLSW}
  J.-P. Cao, H.-Q. Lin, K.-J. Shi and Y. Wang, \textsl{Exact solutions
  and elementary excitations in the XXZ spin chain with unparallel
  boundary fields,} \texttt{cond-mat/0212163}.

\end{thebibliography}
\end{document}